\documentclass[leqno]{conm-p-l} 
\usepackage[bookmarksopen,colorlinks]{hyperref}
\usepackage{url}
\copyrightinfo{2002}{(D. Cox, R. Krasauskas, M. Musta\c t\v a)}

\newtheorem{theorem}{Theorem}[section]

\newtheorem{proposition}[theorem]{Proposition}
\newtheorem{corollary}[theorem]{Corollary}

\newtheorem{definition}[theorem]{Definition}
\newtheorem{example}[theorem]{Example}
\newtheorem{remark}[theorem]{Remark}

\def\CC{{\mathbb C}}
\def\PP{{\mathbb P}}

\def\RR{{\mathbb R}}
\def\ZZ{{\mathbb Z}}

\numberwithin{equation}{section}

\newcommand\res[1]{{\lower1pt\hbox{$|$}}_{\raise.5pt\hbox{${\scriptstyle #1}$}}}

\begin{document}

\title{Universal Rational Parametrizations and Toric Varieties}

\date{\today}

\author[Cox, Krasauskas and Musta\c t\v a]{David Cox, Rimvydas
Krasauskas and Mircea Musta\c t\v a\footnote{This research was
partially done while the third author served as a Clay Mathematics
Institute Long-Term Prize Fellow.}}

\address{Department of Mathematics and Computer Science, Amherst 
College, Amherst, MA 01002, USA}
\email{dac@cs.amherst.edu}

\address{Department of Mathematics and Informatics, Vilnius
University, Naugarduko 14, 2600 Vilnius, Lithuania}
\email{Rimvydas.Krasauskas@maf.vu.lt}

\address{Department of Mathematics, Harvard University, Cambridge, MA
02138, USA}
\email{mustata@math.harvard.edu}

\begin{abstract} This note proves the existence of universal rational
parametrizations.  The description involves homogeneous coordinates on
a toric variety coming from a lattice polytope.  We first describe how
smooth toric varieties lead to universal rational parametrizations of
certain projective varieties.  We give numerous examples and then
discuss what happens in the singular case.  We also describe rational
maps to smooth toric varieties.
\end{abstract}


\maketitle

\section{Introduction}
\label{intro}

In geometric modeling rational curves and surfaces are widely used in
the form of B\'ezier curves and surfaces or simple low-degree
surfaces, e.g., various quadrics, torus surfaces, Dupin cyclides etc.
Construction of curve arcs and patches on a given surface with the
lowest possible parametrization degree is an important task.  For
instance this may help to solve data conversion problems which arise
when translating from traditional solid modeling systems that deal
with such simple surfaces to NURBS-based systems.

It follows that there is a need to understand \emph{all} possible
parametrizations of a given curve or surface.  Is it somehow possible
to find a ``best'' parametrization?  In the case of toric surfaces
(and, more generally, projective toric varieties of any dimension),
this paper will offer one answer to this question, which we call a
\emph{universal rational parametrization}.

To illustrate what we mean by this, we give two examples of surfaces
with universal rational parametrizations.

\begin{example}
\label{firstex}{\rm 
Consider a quadric surface $Q$ given by the homogeneous equation $u_0
u_3 = u_1 u_2$ in projective space $\PP^3$.  Any rational
parametrization of $Q$ can be represented by a collection of
polynomials $H = (h_0,h_{2},h_{2},h_3)$ such that
\begin{equation}
\label{fcond}
h_0 h_3 = h_1 h_2\quad \text{and} \quad
\gcd(h_0,h_{2},h_{2},h_3) = 1.
\end{equation}
One obvious rational parametrization is given by 
\begin{equation}
\label{firstuniv} 
P(x_{1}, x_{2},x_{3},x_{4}) = (x_2 x_3, x_1 x_3, x_2 x_4, x_1 x_4)
\end{equation}
since $(x_2 x_3)(x_1 x_4) = (x_1 x_3)(x_2 x_4)$.

Now suppose that we have a collection of polynomials 
\begin{equation}
\label{hcond}
F = (f_1,f_{2},f_{3},f_4),\quad \gcd(f_1,f_{2}) = \gcd(f_{3},f_4) = 1.
\end{equation}
Then let $H = P\circ F$, i.e.,
\[
H = (h_0,h_{1},h_{2},h_3) = (f_2 f_3, f_1 f_3, f_2 f_4, f_1 f_4).
\]
It is straightforward to show that $H$ satisfies \eqref{fcond}. 
In other words, from the parametrization $P$ of \eqref{firstuniv}, we
get infinitely many others by composing with any $F$ satisfying 
\eqref{hcond}.

But even more is true:\ Theorem~\ref{upthm} implies that \emph{all}
$H$'s satisfying \eqref{fcond} arise in this way.  In other words,
such an $H$ is of the form $H = P\circ F$ for some $F$ as in
\eqref{hcond}.  Furthermore, although $F$ is not unique,
Theorem~\ref{upthm} describes the non-uniqueness precisely:\ given one
$F$ with $H = P\circ F$, then all others are of the form
\[
(\lambda f_1,\lambda f_2,\lambda^{-1} f_3,\lambda^{-1} f_4). 
\]
for some nonzero scalar $\lambda$.}\qed
\end{example}

In the language of Theorem~\ref{upthm}, we say that $P$ from
\eqref{firstuniv} is a \emph{universal rational parametrization} of
the quadric $Q$.  The key property of the quadric $Q$ is that it comes
from $\PP^{1} \times \PP^{1}$.  If $x_{1},x_{2}$ are homogeneous
coordinates on the first $\PP^{1}$ and $x_{3},x_{4}$ are homogeneous
coordinates on the second, then $P$ induces an embedding
\[
\PP^{1} \times \PP^{1} \longrightarrow \PP^{3}
\]
whose image is $Q$.

Here is the second example of a universal rational parametrization.

\begin{example}
\label{firststeiner}{\rm Consider the Steiner surface $S$ in $\PP^3$,
which is defined in homogeneous coordinates by the equation
\[
u_1^2u_2^2+u_2^2u_3^2+u_3^2u_1^2 = u_0u_1u_2u_3.
\]
Note that $S$ is not a smooth surface---its singular locus consists of
the three lines 
\begin{equation}
\label{lines}
u_1 = u_2 = 0,\quad  u_2 = u_3 = 0,\quad  u_3 = u_1 = 0.
\end{equation}
A rational parametrization of $S$ consists of polynomials $H
= (h_0,h_1,h_2,h_3)$ such that
\begin{equation}
\label{stcond1}
h_1^2h_2^2+h_2^2h_3^2+h_3^2h_1^2 = h_0h_1h_2h_3\quad \text{and}\quad
\mathrm{gcd}(h_0,h_1,h_2,h_3) = 1.
\end{equation}
One can easily show that 
\begin{equation}
\label{steinerp}
P(x_1,x_2,x_3) = (x_1^2+x_2^2+x_3^2, x_1x_2 , x_2x_3, x_3x_1)
\end{equation}
is a parametrization of $S$.  Furthermore, given polynomials 
\begin{equation}
\label{stcond}
F = (f_1,f_2,f_3),\quad \gcd(f_1,f_2,f_3) = 1,
\end{equation}
we see that 
\[
H = P \circ F = (f_1^2+f_2^2+f_3^2, f_1f_2 , f_2f_3, f_3f_1)
\]
satisfies \eqref{stcond1} and hence is a rational parametrization of
$S$.   

In this situation, Theorem~\ref{upthm} tells us that \emph{all} $H$'s
satisfying \eqref{stcond1} are of the form $H = P \circ F$ for some
$F$ satisfying \eqref{stcond}, \emph{provided} the image of $H$ does
not lie in the lines \eqref{lines}.  Furthermore, Theorem~\ref{upthm}
implies that $F = (f_1,f_2,f_3)$ is unique up to $\pm 1$.}\qed
\end{example}

By Theorem~\ref{upthm}, \eqref{steinerp} is the universal rational
parametrization of the Steiner surface $S$.  In this case, the key
property of $S$ is that it came from $\PP^2$ via the map $\PP^2 \to S$
induced by $\eqref{steinerp}$.  This map is not an embedding but is
birational (i.e., is generically one-to-one).  Furthermore, the three
lines \eqref{lines} are where the map fails to have an inverse.

Both $\PP^1\times\PP^1$ and $\PP^2$ are examples of \emph{smooth toric
varieties}, and the coordinates $x_1,x_2,x_3,x_4$ for
$\PP^1\times\PP^1$ and $x_1,x_2,x_3$ for $\PP^2$ are examples of
\emph{homogeneous coordinates} of toric varieties.  Hence it makes
sense that there should be a toric generalization of these examples.
For instance, we will see that the gcd conditions \eqref{hcond} and
\eqref{stcond} are dictated by the data which determines the toric
variety.  

The paper is organized into six sections as follows:
\medskip

Section \ref{intro}: Introduction

Section \ref{backgrd}: Background and Related Work

Section \ref{projective}: Universal Rational Parametrizations (Smooth
Case) 

Section \ref{projective}: Universal Rational Parametrizations
(Singular Case)

Section \ref{rational}: Rational Maps to Smooth Toric Varieties

Section \ref{theory}: Theoretical Justifications
\medskip

\noindent In Section \ref{backgrd} we will describe toric varieties
and homogeneous coordinates along with a summary of related work.  In
Section~\ref{projective}, we give a careful definition of rational
parametrization and state Theorems~\ref{upthm}, our main result about
universal rational parametrizations when the toric variety involved is
smooth.  We also give numerous examples.  Then, in
Section~\ref{singular}, we discuss Theorem~\ref{singupthm}, which
describes what happens when the toric variety is singular.  However,
in order to prove these results, we need to understand rational maps
to smooth toric varieties.  This is the subject of
Section~\ref{rational}, where the main result is Theorem~\ref{ratmap}.
Finally, Section~\ref{theory} includes proofs of the results stated in
Sections~\ref{projective}, \ref{singular} and~\ref{rational}.

In this paper, we will work over the complex numbers $\CC$ so that we
can apply the tools of algebraic geometry.  Let $\CC^* =
\CC\setminus\{0\} = \{ z \in \CC \mid z \ne 0\}$.

Geometric modeling is mostly concerned with real varieties.  In
practice, many important real surfaces are real parts of complex toric
surfaces with possibly non-standard real structures.  The results of
this paper hold over $\RR$, provided we use the standard real
structure on the toric varieties involved.  Our results can also be
applied, with some straightforward modifications, to the case of
non-standard real structures.  The details about this situation and
the practical issues of using universal rational parametrizations in
geometric modeling will be presented elsewhere.

We would like to thank the referee for pointing out a problem in our
original version of Theorem~\ref{upthm} and for suggesting the current
form of Example~\ref{hirzebruch}.

\section{Background Material and Related Work}
\label{backgrd}

The concept of universal rational parametrization was introduced at
first for nonsingular quadric surfaces under the name of ``generalized
stereographic projection'' in \cite{DHJ}.  It was extended to more
general rational surfaces in \cite{up} and \cite{M} (see also the
recent paper \cite{Mu}).

Around the same time, homogeneous coordinates for toric varieties
where defined by numerous people---see \cite{hc} for a complete list.
Also important were maps into toric varieties, which were first
explored in \cite{Guest} and \cite{Jac}.  This led the first author to
the description of maps to smooth toric varieties given in
\cite{functor}.

The relation between universal rational parametrizations and toric
varieties was first realized when the second author defined the toric
surface patches in \cite{rimas}.  An account of this may also be found
in \cite{zube}.

\subsection{Toric Varieties}
In this paper, we will assume that the reader is familiar with the
elementary theory of toric varieties, as explained in \cite{what}.  A
toric variety $X_\Sigma$ is determined by a fan $\Sigma$, which is a
collection of cones $\sigma \subset \RR^n$ satisfying certain
properties.  We will assume that the union of the cones in $\Sigma$
is all of $\RR^n$.  This means that $X_\Sigma$ is a compact toric
variety.

Among the cones of $\Sigma$, the edges (= one-dimensional cones) play
a special role.  Suppose that the edges of $\Sigma$ are
$\rho_1,\dots,\rho_r$.  Then each $\rho_i$ corresponds to $x_i$, $n_i$
and $D_i$, where:
\begin{itemize}
\item The variable $x_i$ is in the homogeneous coordinate ring of
$X_\Sigma$. 
\item The vector $n_i \in \ZZ^n$ is the first nonzero integer vector
in $\rho_i$.
\item The subvariety $D_i \subset X_\Sigma$ is defined by $x_i = 0$.
\end{itemize} 
We think of $x_1,\dots,x_r$ as coordinates on $\CC^r$.  We can use
the $x_i$ to construct the toric variety $X_\Sigma$ as follows.  Let
\begin{equation}
\label{gdef}
G = \{(\mu_1,\dots,\mu_r) \in (\CC^*)^r \mid {\textstyle\prod_{i=1}^r}
\mu_i^{\langle m, n_i\rangle} = 1\ \text{for all}\ m \in \ZZ^n\},
\end{equation}
where $\langle \, , \, \rangle$ is dot product on $\ZZ^n$.
This is a subgroup of $(\CC^*)^r$ and hence acts on $\CC^r$ in the
usual way.  Also, for each cone $\sigma \in \Sigma$, let
\[
x^{\hat\sigma} = \prod_{n_i \notin \sigma} x_i
\]
be the product of all variables corresponding to edges \emph{not}
lying in $\sigma$.  Finally, let the \emph{exceptional set} $Z \subset
\CC^r$ be defined by the equations $x^{\hat\sigma} = 0$ for all
$\sigma \in \Sigma$.  Then we get the quotient representation
\begin{equation}
\label{quotient}
X_\Sigma = (\CC^r \setminus Z)/G.
\end{equation}
As explained in \cite{what}, this generalizes the quotient
construction 
\[
\PP^n = (\CC^{n+1} \setminus \{0\})/\CC^*.
\]

One consequence of \eqref{quotient} is that we have a natural map
$\CC^r \setminus Z \to X_\Sigma$.  We can think of this as a
\emph{rational} map from $\CC^r$ to $X_\Sigma$ which is not defined on
the exceptional set $Z$.  We will write this as
\begin{equation}
\label{crxs}
\pi : \CC^r -\! -\!\hskip-1.5pt \to X_\Sigma,
\end{equation}
where the broken arrow means that we have a rational map.  The map
\eqref{crxs} will play an important role in what follows.

\subsection{Polytopes}
In Sections~\ref{projective} and \ref{singular}, we will consider the
projective toric variety $X_\Delta$ determined by an $n$-dimensional
lattice polytope $\Delta \subset \RR^n$.  The idea is that for each
face of $\Delta$, we get the cone generated by the inward-pointing
normals of the facets of $\Delta$ containing the face.  This gives the
\emph{normal fan} $\Sigma_\Delta$ of $\Delta$, and the corresponding
toric variety $X_{\Sigma_\Delta}$ is denoted $X_\Delta$.

Observe that edges of the normal fan correspond to facets of $\Delta$.
Hence the homogeneous coordinates $x_1,\dots,x_r$ correspond to the
facets of $\Delta$.  For this reason, we call $x_1,\dots,x_r$ the
\emph{facet variables} of the polytope $\Delta$.

We can use $\Delta$ to obtain some interesting monomials in the facet
variables.  Represent $\Delta$ as the intersection
\begin{equation}
\label{deltadesc}
\Delta = \bigcap_{i=1}^{r} \{m \in \RR^n \mid 
\langle m, n_{i}\rangle \ge -a_{i}\}
\end{equation}
of closed half-spaces.  This gives the following monomials and
polynomials.

\begin{definition}
\label{dmsd}
For each lattice point $m \in \Delta\cap
\ZZ^n$, we define the {\bfseries $\Delta$-monomial} $x^{m}$ to be
\begin{equation}
\label{deltam}
x^{m} = \prod_{i=1}^{r} x_{i}^{\langle m, n_{i}\rangle+a_{i}}.
\end{equation}
We also define $S_\Delta$ to be the linear span of the set of
$\Delta$-monomials.  Thus
\[
S_\Delta = \mathrm{Span}(x^{m} \mid m \in \Delta\cap \ZZ^n).
\]
\end{definition}

Since the $i$th facet is defined by $\langle m, n_{i}\rangle+a_{i} =
0$ and $\langle m, n_{i}\rangle+a_{i} \ge 0$ on $\Delta$ ($n_i$ points
inward), we see that the exponent of $x_i$ measures the ``distance''
(in the lattice sense) from $m$ to the $i$th facet.

Here is an example of facet variables and $\Delta$-monomials.

\begin{example}
\label{p1p1b2first}
{\rm Consider the polytope $\Delta$
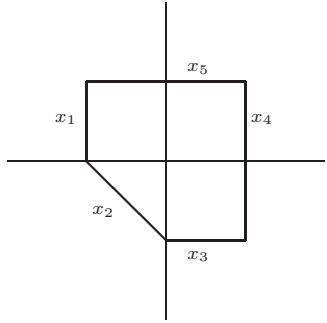
\[
\begin{matrix}
\begin{picture}(120,120)
\put(0,60){\line(1,0){120}}
\put(60,0){\line(0,1){120}}
\thicklines
\put(30,60){\line(0,1){30}}
\put(30,90){\line(1,0){60}}
\put(90,90){\line(0,-1){60}}
\put(90,30){\line(-1,0){30}}
\put(60,30){\line(-1,1){30}}
\put(18,75){$\scriptstyle{x_1}$}
\put(68,94){$\scriptstyle{x_5}$}
\put(32,40){$\scriptstyle{x_2}$}
\put(92,75){$\scriptstyle{x_4}$}
\put(68,23){$\scriptstyle{x_3}$}
\end{picture}
\end{matrix}
\]
with vertices $(1,1), (-1,1), (-1,0), (0, -1), (1,-1)$.  In terms of
\eqref{deltadesc}, we have $a_1 = \dots = a_5 = 1$.  This gives a
toric surface $X_\Delta$ with variables $x_1,\dots,x_5$ as indicated
in the picture.  The 8 points $m \in \Delta\cap\ZZ^2$ give the
following 8 $\Delta$-monomials $x^m$:
\[
\begin{array}{cccccl}

 x_2 x_3^2 x_4^2, & x_1 x_2^2 x_3^2 x_4, & x_1^2 x_2^3 x_3^2, \\
 x_3 x_4^2 x_5, & x_1 x_2 x_3 x_4 x_5, & x_1^2 x_2^2 x_3 x_5, \\
 & x_1 x_4 x_5^2, & x_1^2 x_2 x_5^2 .
\end{array}
\]
We will say more about this example in Sections~\ref{projective}
and~\ref{rational}.}\qed 
\end{example}

We should also mention that polynomials $q \in S_\Delta$ have the
following important property:\ given $\mu = (\mu_1,\dots,\mu_r)$ in
the group $G$ defined in \eqref{gdef}, one easily sees that
\begin{equation}
\label{pmx}
q(\mu_1 x_1,\dots,\mu_r x_r) = \mu_\Delta\, q(x_1,\dots,x_r),
\end{equation}
where 
\begin{equation}
\label{mud}
\mu_\Delta = \prod_{i=1}^r \mu_i^{a_i}.
\end{equation}

\subsection{Rational Maps to Projective Space}
Now pick a collection $P = (p_0,\dots,p_s)$ of $s+1$ polynomials in
$S_\Delta$.  This gives a rational map
\[
p: \CC^r -\! -\!\hskip-1.5pt \to \PP^s.
\]
If $X$ denotes the Zariski closure of the image, then we write $p$ as
\[
p: \CC^r -\! -\!\hskip-1.5pt \to X \subset \PP^{s}.
\]
We can relate $p$ to the toric variety $X_{\Delta}$ as follows.

\begin{proposition}
\label{pfactors}
In the above situation, the map $p$ factors $p = \Pi\circ \pi$, 
where $\pi : \CC^{r} -\! -\!\hskip-1.5pt \to X_{\Delta}$ is from 
\eqref{crxs} and 
\[
\Pi: X_\Delta  -\! -\!\hskip-1.5pt \to X
\]
is a rational map.
\end{proposition}

\begin{proof} Given $a = (a_1,\dots,a_r) \in \CC^r\setminus
Z$ and $\mu = (\mu_1,\dots,\mu_r) \in G$, then we get 
$\mu\cdot a = (\mu_1a_1,\dots,\mu_ra_r)$.  By \eqref{pmx}, we have
\[
(p_0(\mu\cdot a),\dots,p_s(\mu\cdot a)) = \mu_\Delta
(p_0(a),\dots,p_s(a)). 
\]
This shows that $p$ induces $\Pi: (\CC^r\setminus Z)/G -\!
-\!\hskip-1.5pt \to \PP^s$.  By \eqref{quotient}, we can identify the
quotient with $X_{\Sigma}$, and the proposition follows.
\end{proof}

When $X_\Delta$ is smooth and $\Pi : X_\Delta -\!  -\!\hskip-1.5pt \to
X$ is sufficiently nice, Theorem~\ref{upthm} asserts that $p$ is a
universal rational parametrization.  In Section~\ref{projective}, we
will use this theorem to explain Examples~\ref{firstex}
and~\ref{firststeiner}.

We will also see in Section~\ref{singular} that this doesn't quite
work when $X_\Delta$ is singular.  In this case,
Theorem~\ref{singupthm} will show that we get a universal rational
parametrization by considering a suitable resolution of singularities.

\section{Universal Rational Parametrizations (Smooth Case)}
\label{projective}

In this section, we will prove the existence of universal rational
parametrizations for certain projective varieties which arise
naturally from smooth toric varieties associated to polytopes.  

\subsection{Rational Parametrizations}

We first give a definition of rational pa\-rametrization which is
useful in geometric modeling.  Given a projective variety $Y \subset
\PP^{k}$, we define its \emph{affine cone} $C_{Y} \subset \CC^{k+1}$
to be
\[
C_{Y} = \pi^{-1}(Y) \cup \{0\} \subset \CC^{k+1},
\]
where $\pi : \CC^{k+1} \setminus \{0\} \to \PP^{k}$ is the usual map.
Using this, we can make the following definition.

\begin{definition}
\label{rpdef}
Let $R = \CC[y_1,\dots,y_d]$ be the coordinate ring of $\CC^d$.
A {\bfseries rational parametrization} of a projective variety $Y 
\subset \PP^{s}$ consists of $H = (h_{0},\dots,h_{s}) \in R^{s+1}$ 
such that $\mathrm{gcd}(h_{0},\dots,h_{s}) = 1$ and $H(\CC^{d}) 
\subset C_{Y}$.
\end{definition}

In this paper, we use the convention that $\mathrm{gcd}(0,\dots,0) =
0$.  Hence the gcd condition implies that the polynomials in a
rational parametrization are not all zero.  Then $H(\CC^{d}) \subset
C_{Y}$ implies that $H : \CC^d \to \CC^{s+1}$ induces a rational map
\[
h : \CC^d  -\! -\!\hskip-1.5pt \to \PP^{s}
\]
whose image lies in $Y$.  It is important to note that in
Definition~\ref{rpdef}, we do not require that $h : \CC^d -\!
-\!\hskip-1.5pt \to Y$ be surjective or have dense image.  Thus a
rational parametrization might only parametrize a proper subvariety of
$Y$.  Also note that the gcd condition of Definition~\ref{rpdef}
implies that two rational parametrizations $H$ and $H'$ give the same
rational map to $\PP^{s}$ if and only if $H = cH'$ for $c \ne 0$ in
$\CC$.

\subsection{One Particular Parametrization}  Let $\Delta$ be an
$n$-dimensional lattice polytope in $\RR^n$.  This gives the 
toric variety $X_\Delta$ determined by the normal fan $\Sigma_\Delta$
of $\Delta$.  We will assume that $X_\Delta$ is smooth.

By Definition~\ref{dmsd}, the facet variables 
$x_{1},\dots,x_{r}$ and the lattice points $\Delta\cap
\ZZ^{n}$ give rise to the vector space of polynomials
\[
S_{\Delta} = \mathrm{Span}(x^{m} \mid m \in \Delta \cap \ZZ^{n}),
\]
where $x^{m}$ is the $\Delta$-monomial.  As in Section~\ref{backgrd},
a collection of $s+1$ polynomials
\begin{equation}
\label{firstp}
P = (p_{0},\dots,p_{s}) \in S_{\Delta}^{s+1}
\end{equation}
gives a rational map
\[
p : \CC^{r} -\!-\!\hskip-1.5pt \to X \subset \PP^{s}
\]
where $X$ is the Zariski closure of the image.  Then
Proposition~\ref{pfactors} shows that $p = \Pi\circ\pi$, where $\pi :
\CC^{r} -\!-\!\hskip-1.5pt \to X_{\Delta}$ is from \eqref{crxs} and
\[
\Pi : X_{\Delta} -\!-\!\hskip-1.5pt \to X \subset \PP^s
\]
is a rational map.  As already mentioned, the basic idea is that $P$
is a universal rational parametrization when $\Pi$ is sufficiently
nice.  However, we need to make some further definitions before we can
state our main result.

\subsection{Sufficiently Nice} We can finally explain when
$\Pi : X_\Delta -\!-\!\hskip-1.5pt \to X$ is sufficiently nice.  Using
the above notation, this means the following two things:
\begin{itemize}
\item First, $\Pi$ is \emph{strictly defined}, which means for every
  $a \in \CC^r \setminus Z$, we have $p_i(a) \ne 0$ for some $0 \le i
  \le s$.  Using $X_\Delta = (\CC^d \setminus Z)/G$ and
  Proposition~\ref{pfactors}, one can show that this condition ensures
  that $\Pi$ is defined everywhere.  Thus we write $\Pi : X_\Delta \to
  X$ when $\Pi$ is strictly defined.
\item Second, $\Pi$ is \emph{birational}, which means that
  $\Pi$ induces an isomorphism between dense open subsets of
  $X_\Delta$ and $X$.
\end{itemize}
When we discuss projections later in the section, we will give several
conditions which are equivalent to being strictly defined.  However,
we will see in Example~\ref{hirzebruch} that being strictly defined is
in general \emph{stronger} than just assuming that the rational map
$\Pi$ is defined everywhere on $X_\Delta$.

An important observation is that the $p_{i}$ in \eqref{firstp} are
relatively prime when $\Pi$ is strictly defined.  To prove this,
suppose that some nonconstant polynomial $q$ divides the $p_{i}$.
Since the exceptional set $Z \subset \CC^{r}$ has codimension at least
$2$, we can find $a \in \CC^{r} \setminus Z$ such that $q(a) = 0$.
Hence $p_{i}(a) = 0$ for all $i$, which can't happen when $\Pi$ is
strictly defined.  This proves that the $p_{i}$ are relatively prime.
By Definition~\ref{rpdef}, it follows that $P$ is a rational
parametrization of $X$.

We also note that being strictly defined implies that $\Pi$ and
hence $p$ are onto, i.e., $X$ is the image of $p : \CC^{r}
-\!-\!\hskip-1.5pt \to \PP^{s}$.  This follows because $X_\Delta$ is
compact.

Finally, note that when $\Pi$ is birational, there is a nonempty
Zariski open subset
\[
U \subset X
\]
on which $\Pi^{-1}$ is defined.  We may assume that $U$ is the maximal
such open set.

\subsection{$\Sigma_\Delta$-Irreducible Polynomials} The rough idea of
a universal rational parametrization is that any rational
parametrization $H$ should arise from $P$ by composition with a
polynomial map $\CC^d \to \CC^r$.  But if the image of $\CC^d \to
\CC^r$ lies in the exceptional set $Z$, then the composition doesn't
make sense since $p$ is not defined on $Z \subset \CC^r$.  It follows
that we need to exclude certain polynomial maps.  The precise
definition is as follows.  Let $R = \CC[y_1,\dots,y_d]$.  

\begin{definition}
\label{irreddef}
We say that $F = (f_1,\dots,f_r) \in R^r$ is {\bfseries
$\Sigma_\Delta$-irreducible} if $\mathrm{gcd}(f_{i_1},\dots,f_{i_k}) =
1$ whenever no cone of $\Sigma_\Delta$ contains
$\rho_{i_1},\dots,\rho_{i_k}$.
\end{definition}

Because of our gcd convention, Definition~\ref{irreddef} implies in
particular that the edges $\rho_i$ such that $f_i = 0$ all lie in some
cone of $\Sigma_\Delta$.  In the discussion which follows, we will
identify $F$ with the polynomial function $\CC^d \to \CC^r$ it
induces.

Here are two examples of this definition.

\begin{example}
\label{p2irred}{\rm Consider the toric variety $X_\Delta = \PP^2$
coming from the polytope $\Delta$ with vertices $(0,0), (2,0), (0,2)$.
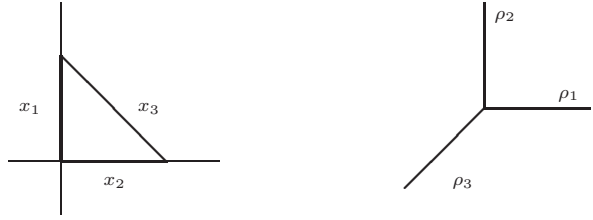
\[
\begin{matrix}
\begin{picture}(260,80)
\thinlines
\put(20,20){\line(1,0){80}}
\put(40,0){\line(0,1){80}}
\thicklines
\put(40,20){\line(1,0){40}}
\put(40,20){\line(0,1){40}}
\put(40,60){\line(1,-1){40}}
\put(69,39){$\scriptstyle{x_3}$}
\put(24,39){$\scriptstyle{x_1}$}
\put(56,11){$\scriptstyle{x_2}$}
\put(200,40){\line(1,0){40}}
\put(200,40){\line(0,1){40}}
\put(200,40){\line(-1,-1){30}}
\put(228,44){$\scriptstyle{\rho_1}$}
\put(204,74){$\scriptstyle{\rho_2}$}
\put(188,11){$\scriptstyle{\rho_3}$}
\end{picture}
\end{matrix}
\] 
The polytope $\Delta$ is on the left with facet variables
$x_1,x_2,x_3$ and the normal fan is on the right with edges
$\rho_1,\rho_2,\rho_3$.  The only choice for
$\rho_{i_1},\dots,\rho_{i_k}$ in Definition~\ref{irreddef} is
$\rho_1,\rho_2,\rho_3$, so that $F = (f_1,f_2,f_3)$ is
$\Sigma$-irreducible if $\gcd(f_1,f_2,f_3) = 1$.  This is the gcd
condition \eqref{stcond} in Example~\ref{firststeiner}. }\qed
\end{example}

\begin{example}
\label{p1p1irred}{\rm Let $\Delta$ be the unit square in the plane
with vertices $(0,0)$, $(1,0)$, $(1,1)$, $(0,1)$.   This gives 
$X_\Delta = \PP^1\times\PP^1$.
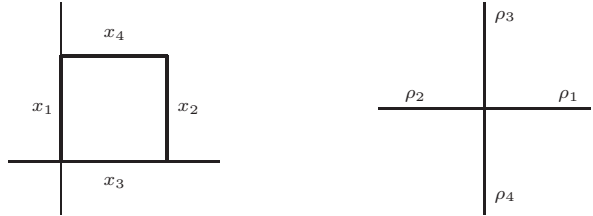
\[
\begin{matrix}
\begin{picture}(260,80)
\thinlines
\put(20,20){\line(1,0){80}}
\put(40,0){\line(0,1){80}}
\thicklines
\put(40,20){\line(1,0){40}}
\put(40,20){\line(0,1){40}}
\put(40,60){\line(1,0){40}}
\put(80,20){\line(0,1){40}}
\put(84,39){$\scriptstyle{x_2}$}
\put(56,67){$\scriptstyle{x_4}$}
\put(29,39){$\scriptstyle{x_1}$}
\put(56,11){$\scriptstyle{x_3}$}
\put(160,40){\line(1,0){80}}
\put(200,0){\line(0,1){80}}
\put(228,44){$\scriptstyle{\rho_1}$}
\put(204,74){$\scriptstyle{\rho_3}$}
\put(170,44){$\scriptstyle{\rho_2}$}
\put(204,6){$\scriptstyle{\rho_4}$}
\end{picture}
\end{matrix}
\]
As in the previous example, $\Delta$ and the facet variables
$x_1,x_2,x_3,x_4$ are on the left and $\Sigma_\Delta$ and the edges
are $\rho_1,\rho_2,\rho_3,\rho_4$ are on the right.  The minimal
choices for $\rho_{i_1},\dots,\rho_{i_k}$ in Definition~\ref{irreddef}
are $\rho_1,\rho_2$ and $\rho_3,\rho_4$ (you should check this).  Thus
$F = (f_1,f_2,f_3,f_4)$ is $\Sigma$-irreducible if $\gcd(f_1,f_2) =
\gcd(f_3,f_4) = 1$.  This is the gcd condition \eqref{hcond} in
Example~\ref{firstex}. }\qed
\end{example}

\subsection{The Main Result}
Before stating our main result, we need to introduce some notation.
As above, let $R = \CC[y_{1},\dots,y_{d}]$ be the coordinate ring of
$\CC^{d}$.  Also, given polynomials $F = (f_{1},\dots,f_{r}) \in R^r$
and $m \in \Delta\cap \ZZ^n$, we set
\[
f^{m} = \prod_{i=1}^{r} f_{i}^{\langle m,n_{i}\rangle+a_{i}}.
\]
Recall the group $G$ from \eqref{gdef} and that $\mu \in G$ gives
$\mu_\Delta \in \CC^*$ defined in \eqref{mud}.  The map $\mu \mapsto
\mu_\Delta$ is a group homomorphism $G \to \CC^*$.  Let $G_\Delta$ be
the kernel of this map.  This group will measure the lack of
uniqueness in Theorem~\ref{upthm}.  Finally, let $\sum_m$ denote
summation over all $m \in \Delta\cap \ZZ^n$.

Here is our precise result.

\begin{theorem}
\label{upthm}
Let ${P} = (p_0,\dots,p_s) = (\sum_{m} a_{0m} x^{m}, \dots,\sum_{m}
a_{sm}x^{m}) \in S_\Delta^{s+1}$, where $X_\Delta$ is smooth and $\Pi
: X_\Delta \to X$ is strictly defined and birational, with an inverse
defined on $U \subset X$ which we assume to be maximal.  Then $P$ is
a {\bfseries universal rational parametrization} of $X \subset \PP^s$
in the following sense:
\begin{enumerate}
\item If $F = (f_{1},\dots,f_{r}) \in R^{r}$ is
$\Sigma_{\Delta}$-irreducible, then
\[
{P}\circ F = ({\textstyle \sum_{m} a_{0m} f^{m}, \dots,\sum_{m} 
a_{sm}f^{m}}) \in R^{s+1}
\]
is a rational parametrization of $X \subset \PP^{s}$.
\item Conversely, given any rational parametrization $H \in R^{s+1}$
of $X$ whose image meets $U \subset X$, there is a
$\Sigma_{\Delta}$-irreducible $F = (f_{1},\dots,f_{r}) \in R^{r}$ such
that $H = {P}\circ F$.
\item If $F$ and $F'$ are $\Sigma_{\Delta}$-irreducible, then $P\circ
F = P\circ F'$ as rational parametrizations if and only if $F' =
\mu\cdot F$ for some $\mu \in G_\Delta$.
\end{enumerate}
\end{theorem}

The proof will be given in Section \ref{theory}.  Here is a corollary
of Theorem~\ref{upthm}.

\begin{corollary}
\label{univcor}
Assume the same hypothesis as Theorem~\ref{upthm} and suppose that
$H' = (h_0',\dots,h_s') \in R^{s+1}$ gives a rational map $\CC^d -\!
-\!\hskip-1.5pt \to \PP^{s}$ whose image lies in $X$ and meets $U$.
Then there is a polynomial $q \in R$ and a $\Sigma_\Delta$-irreducible
$F = (f_1,\dots,f_r) \in R^r$ such that
\[
H' = q\, P \circ F.
\]
\end{corollary}

\begin{proof}
Write $H' = q\, H$, where the entries of $H$ are relatively prime.
Since $H$ is a rational parametrization, we are done by
Theorem~\ref{upthm}.
\end{proof}

\subsection{Embeddings}

In order for Theorem~\ref{upthm} to be useful, we need to have a good 
supply of parametrizations $P = (p_0,\dots,p_s) \in S_\Delta^{s+1}$
which satisfy the hypotheses of the theorem.  The first crucial
observation is that since $X_\Delta$ is a smooth toric variety, it is
a standard result that the collection of \emph{all} $\Delta$-monomials
gives a projective embedding (see \cite[Sec.\ 3.4]{Fulton}).  

This means the following.  Suppose that $\Delta\cap\ZZ^n =
\{m_0,\dots,m_\ell\}$ and let
\begin{equation}
\label{pddef}
P_\Delta = (x^{m_0},\dots,x^{m_\ell})
\end{equation}
Then $P_\Delta$ induces $p_\Delta : \CC^r -\!-\!\hskip-1.5pt \to
\PP^{\ell}$, and in the factorization $p_\Delta = \pi\circ \Pi_\Delta$ of
Proposition~\ref{pfactors}, the map $\Pi : X_\Delta \to \PP^{\ell}$ is
an embedding.  Hence we can write $X_\Delta \subset \PP^\ell$.

All of the hypotheses of Theorem~\ref{upthm} are satisfied in this
situation, and the open set $U$ is all of $X = X_\Delta$.  Thus
$P_\Delta$ is a universal rational parametrization in the strong sense
that \emph{every} rational parametrization is of the form
$P_\Delta\circ F$ for some $\Sigma_\Delta$-irreducible $F$.

Here are two examples of this result.

\begin{example}
\label{p1p1second}{\rm Let $\Delta$ be the unit square from
Example~\ref{p1p1irred}.  
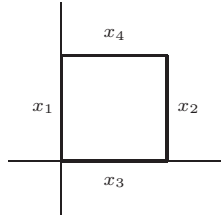
\[
\begin{matrix}
\begin{picture}(120,80)
\put(20,20){\line(1,0){80}}
\put(40,0){\line(0,1){80}}
\thicklines
\put(40,20){\line(1,0){40}}
\put(40,20){\line(0,1){40}}
\put(40,60){\line(1,0){40}}
\put(80,20){\line(0,1){40}}
\put(84,39){$\scriptstyle{x_2}$}
\put(56,67){$\scriptstyle{x_4}$}
\put(29,39){$\scriptstyle{x_1}$}
\put(56,11){$\scriptstyle{x_3}$}
\end{picture}
\end{matrix}
\]
The labeling of $x_1,x_2,x_3,x_4$ is consistent with
Example~\ref{p1p1irred}.  In terms of \eqref{deltadesc}, $a_1 = a_3 =
0$ and $a_2 = a_4 = 1$.  Then
\[
\begin{array}{ccccl}
P_\Delta &=& 
(x_2 x_3, & x_1 x_3, & \\
 &&\, x_2 x_4, & x_1 x_4)& \in S_\Delta^4
\end{array}
\]
gives a universal rational parametrization of its image in $\PP^3$,
which is the quadric $Q$ of Example~\ref{firstex}.  This means that
any parametrization of $Q$ is of the form $P_\Delta\circ F$, where $F
= (f_1,f_2,f_3,f_4)$ satisfies the gcd condition worked out in
Example~\ref{p1p1irred}.

To study uniqueness, we need to compute
\[
G = \{(\mu_1,\mu_2,\mu_3,\mu_4) \in (\CC^*)^4 \mid
{\textstyle\prod_{i=1}^4} 
\mu_i^{\langle m, n_{i}\rangle} = 1\ \text{for all}\ m \in \ZZ^4\}.
\]
Since it suffices to use $m = (1,0)$ and $(0,1)$, we see that 
\[
G = \{(\mu_1,\mu_2,\mu_3,\mu_4) \in (\CC^*)^4 \mid 
\mu_{1}\mu_{2}^{-1} = \mu_{3}\mu_{4}^{-1} =1\},
\]
and it follows that
\[
G =\{(\mu_1,\mu_1,\mu_2,\mu_2) \mid \mu_1, \mu_2 \in \CC^*\}.
\]
Then $a_1 = a_3 = 0$ and $a_2 = a_4 = 1$ imply that if $\mu =
(\mu_1,\mu_1,\mu_2,\mu_2) \in G$, then $\mu_\Delta = \mu_1\mu_2$.
Hence 
\[
G_\Delta = \{(\lambda,\lambda,\lambda^{-1},\lambda^{-1}) \mid \lambda
\in \CC^*\}.  
\]
Thus, when we write a rational parametrization of $Q$ as
$P_\Delta\circ F$, we see that $F = (f_1,f_2,f_3,f_4)$ is unique up to
\[
(\lambda f_1,\lambda f_2,\lambda^{-1} f_3,\lambda^{-1} f_4). 
\]
for some nonzero scalar $\lambda$.  It follows that we obtain
precisely the description given in Example~\ref{firstex}.}\qed
\end{example}

Here is how Theorem~\ref{upthm} applies to one of our earlier
examples.

\begin{example}
\label{p1p1b2}
{\rm Consider the toric variety $X_\Delta$ of the polytope $\Delta$
with vertices $(1,1), (-1,1), (-1,0), (0, -1), (1,-1)$ from
Example~\ref{p1p1b2first}. 
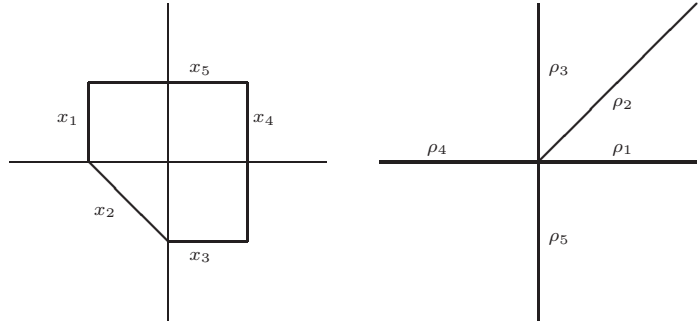
\[
\begin{matrix}
\begin{picture}(260,120)
\put(0,60){\line(1,0){120}}
\put(60,0){\line(0,1){120}}
\thicklines
\put(30,60){\line(0,1){30}}
\put(30,90){\line(1,0){60}}
\put(90,90){\line(0,-1){60}}
\put(90,30){\line(-1,0){30}}
\put(60,30){\line(-1,1){30}}
\put(18,75){$\scriptstyle{x_1}$}
\put(68,94){$\scriptstyle{x_5}$}
\put(32,40){$\scriptstyle{x_2}$}
\put(92,75){$\scriptstyle{x_4}$}
\put(68,23){$\scriptstyle{x_3}$}
\put(140,60){\line(1,0){120}}
\put(200,0){\line(0,1){120}}
\put(200,60){\line(1,1){60}}
\put(228,64){$\scriptstyle{\rho_1}$}
\put(204,94){$\scriptstyle{\rho_3}$}
\put(228,81){$\scriptstyle{\rho_2}$}
\put(158,64){$\scriptstyle{\rho_4}$}
\put(204,30){$\scriptstyle{\rho_5}$}
\end{picture}
\end{matrix}
\]
As usual, the polytope is on the left and the normal fan of $\Delta$
is on the right.  One can show that $X_\Delta$ is the blowup of $\PP^1
\times \PP^1$ at one point.

In terms of \eqref{deltadesc} and the above labeling, we have $a_1 =
\dots = a_5 = 1$.  The 8 points of $\Delta\cap \ZZ^2$ give an
embedding of $X_\Delta$ into $\PP^7$.  It follows that
\[
\begin{array}{cccccl}
P_\Delta &=& 
(x_2 x_3^2 x_4^2, & x_1 x_2^2 x_3^2 x_4, & x_1^2 x_2^3 x_3^2,& \\
 && x_3 x_4^2 x_5, & x_1 x_2 x_3 x_4 x_5, & x_1^2 x_2^2 x_3 x_5,& \\
 &&& x_1 x_4 x_5^2, & x_1^2 x_2 x_5^2)& \in S_\Delta^8.
\end{array}
\]
is the universal rational parametrization of $X_\Delta$ by
Theorem~\ref{upthm}.  

Let's work out what this means.  According to
Definition~\ref{irreddef}, $F = (f_{1},\dots,f_{5})$ is
$\Sigma_\Delta$-irreducible provided that
\[
\gcd(f_{1},f_{3}) = \gcd(f_{1},f_{4}) = \gcd(f_{2},f_{4}) = 
\gcd(f_{2},f_{5}) = \gcd(f_{3},f_{5}) = 1.
\]
Then any rational parametrization of $X_{\Sigma}$ is of the form
$P_\Delta\circ F$ for some $F$ satisfying this gcd condition.  To
determine the lack of uniqueness, we need to compute the group $G$.
Using the methods of Example~\ref{p1p1second}, one obtains
\[
G = \{(\lambda,\mu,\nu,\lambda \mu, \mu\nu) \mid \lambda,\mu,\nu \in 
\CC^{*}\},
\]
and then $a_1 = \dots = a_5 = 1$ imply that
\[
G_\Delta = \{(\lambda,\mu,\nu,\lambda \mu, \mu\nu) \mid \lambda,\mu,\nu \in 
\CC^{*}, \lambda^2 \mu^3 \nu^2 = 1\}.
\]
Thus rational parametrizations of $X_{\Sigma}$ are all of the form
$P_\Delta\circ F$, where $F$ is unique up to $(\lambda,\mu,\nu,\lambda \mu,
\mu\nu)\cdot F$ for $\lambda^2 \mu^3 \nu^2 = 1$.\qed}
\end{example}

\subsection{Projections}

Although $P_\Delta = (x^{m_0},\dots,x^{m_\ell})$ from \eqref{pddef}
always gives a universal rational parametrization, it is rarely useful
in practice since it usually gives an embedding into a projective
space of high dimension.  An important observation is that we can
think of the general case
\[
P = (p_0,\dots,p_s) = ({\textstyle\sum_{i=0}^\ell a_{0i} x^{m_i},
\dots,\sum_{i=0}^\ell a_{si}x^{m_i}}) \in S_\Delta^{s+1}
\]
in terms of projections.  Let $\PP^\ell$ be a projective space with
homogeneous coordinates $z_0,\dots,z_\ell$.  Then the $s+1$ linear
forms $\sum_{i=0}^\ell a_{ji} z_i$ define a projection $\PP^{\ell}
-\!-\!\hskip-1.5pt \to \PP^{s}$ defined by
\[
(z_0,\dots,z_\ell) \mapsto ({\textstyle\sum_{i=0}^\ell a_{0i} {z_i},
\dots,\sum_{i=0}^\ell a_{si}{z_i}}).
\]
The \emph{center} $L \subset \PP^{\ell}$ of this projection is
defined by $\sum_{i=0}^\ell a_{ji} z_i = 0$ for $j = 0,\dots,s$.  This
tells us where the projection is not defined.

If we compose this projection with $p_\Delta : \CC^r -\!
-\!\hskip-1.5pt \to \PP^\ell$, then we get the rational map $p : \CC^r
-\!  -\!\hskip-1.5pt \to \PP^s$ induced by $P$.  Furthermore, since
the image of $p_\Delta$ is $X_\Delta \subset \PP^\ell$, it follows
that the variety $X$ parametrized by $P$ is the image of $X_\Delta$
under the projection.

From this point of view, we can think of $\Pi$ as a projection.  It is
then straightforward to check that $\Pi$ is strictly defined if and
only if $X_{\Delta}$ is disjoint from the center $L$ of the
projection.  (For more sophisticated readers, we point out that being
strictly defined is equivalent to the assertion that the linear system
on $X_\Delta$ spanned by the $p_i$ has no base points.  One can also
show that $X_\Delta$ is the normalization of $X$ when $\Pi$ is
strictly defined and birational.)

Let's give an example from geometric modeling which involves the 
projection of a toric variety.

\begin{example}{\rm Consider the toric variety $X_\Delta = \PP^2$,
where $\Delta$ is the polytope from Example~\ref{p2irred}.  In terms
of \eqref{deltadesc}, we have $a_1 = a_2 = 0$ and $a_3 = 2$.  The 6
points of $\Delta \cap \ZZ^2$ define
\[
P_\Delta = (x_1^2, x_2^2, x_3^2, x_1 x_2, x_2 x_3, x_3 x_1).
\]
This gives the usual Veronese embedding of $\PP^2$ into $\PP^5$.

The composition of this map with the projection $\PP^5
-\!-\!\hskip-1.5pt \to \PP^3$ defined by
\begin{equation}
\label{stpr}
(z_0,z_1,z_2,z_3,z_4,z_5) \mapsto (z_0+z_1+z_2,z_3,z_4,z_5) 
\end{equation}
gives a rational parametrization
\begin{equation}
\label{upsteiner}
P = (x_1^2+x_2^2+x_3^2, x_1x_2 , x_2x_3, x_3x_1)
\end{equation}
of the Steiner surface $S \subset \PP^3$ defined by
\[
u_1^2u_2^2+u_2^2u_3^2+u_3^2u_1^2 = u_0u_1u_2u_3,
\]
where $u_0,u_1,u_2,u_3$ are homogeneous coordinates on $\PP^3$.  We
saw this equation earlier in Example~\ref{firststeiner}.

It is easy to check that the center of the projection \eqref{stpr} is
disjoint from $X_\Delta$.  Thus the map $\Pi : X_\Delta \to S$ is
strictly defined.  Furthermore, since
\[
x_1^2 = \frac{(x_1x_2)(x_3x_1)}{x_2x_3} = \frac{u_1u_3}{u_2},
\]
one easily sees that $\Pi : X_\Delta \to S$ is birational with inverse 
\begin{align*}
\Pi^{-1}(u_0,u_1,u_2,u_3) &= \Big(\frac{u_1u_3}{u_2}, 
\frac{u_1u_2}{u_3},\frac{u_2u_3}{u_1},u_1,u_2,u_3\Big)\\
	&= (u_1^2 u_3^2, u_1^2 u_2^2, u_2^2u_3^2, u_1^2u_2u_3,
u_1u_2^2u_3, u_1u_2u_3^2).
\end{align*}
Also notice that $\Pi^{-1}$ is defined on the complement of the three
lines $u_1 = u_2 = 0$, $u_2 = u_3 = 0$, $u_3 = u_1 = 0$.

By Theorem~\ref{upthm}, \eqref{upsteiner} is the universal rational
parametrization of the Steiner surface $S$.  It follows that if $H$ is
a rational parametrization of $S$ whose image is not contained in any
of the above three lines, then $H = P\circ F$, where $F =
(f_1,f_2,f_3)$.  Furthermore, we know that $F$ is
$\Sigma_\Delta$-irreducible, which by Example~\ref{p2irred} means
$\gcd(f_1,f_2,f_3) = 1$.

Finally, we know that $G = \{(\lambda,\lambda,\lambda) \mid \lambda
\in \CC^*\} \simeq \CC^*$ in this case.  Since $a_1 = a_2 = 0$ and
$a_3 = 2$, we see that $\mu_\Delta = \lambda^2$ when $\mu =
(\lambda,\lambda,\lambda)$.  It follows that the kernel of $\mu
\mapsto \mu_\Delta$ is $\pm(1,1,1)$, so that in $H = P\circ F$, $F$ is
unique up to $\pm1$.  Hence we recover the description of the rational
parametrizations of the Steiner surface given in
Example~\ref{firststeiner}.

Observe that $F = (f_1,f_2,f_3)$ may fail to exist when the image of
$H$ is contained in one of the three lines $z_1 = z_2 = 0$, $ z_2 =
z_3 = 0$, $z_3 = z_1 = 0$.  For example, $H = (u,0,0,v)$ is a rational
parametrization from $\CC^2$ to $S$ which is not of the form $P\circ
F$ for any $F \in \CC[u,v]^3$.  Notice also that the union of these
lines is the singular locus of $S$.}\qed
\end{example}

We next describe one important class of projections which always lead
to universal rational parametrizations.  Suppose that the smooth
$n$-dimensional toric variety $X_\Delta$ is embedded into $\PP^\ell$
via $P_\Delta$.  Then let
\[
P = (p_0,\dots,p_{n+1}) \in S_\Delta^{n+2}
\]
be chosen generically.  Then $X \subset \PP^{n+1}$ is the image of
$X_\Delta$ under a generic projection.  It is well-known that in this
situation, $X_\Delta$ is disjoint from the center of the projection
and the restriction of the projection to $X_\Delta$ is birational.
Hence $P$ is a universal rational parametrization in this generic
case.  Notice that $X$ has dimension $n$ and hence is a hypersurface
in $\PP^{n+1}$.

In particular, when  $X_\Delta$ is a smooth toric surface, it follows
that
\[
P = (p_0,\dots,p_3) \in S_\Delta^{4}
\]
is a universal rational parametrization whenever the $p_i$ are chosen
generically.  Here, we parametrize a surface in $\PP^3$, which is the
case of greatest interest in geometric modeling.

However, we should also mention that there are some non-generic
projections which also work nicely.  Here is another class of
projections which are guaranteed to give universal rational
parametrizations.

\begin{proposition}
\label{lattice}
Let $X_\Delta$ be the smooth toric variety of a polytope $\Delta$ and
let $\mathcal{A} = \{\tilde{m}_0,\dots,\tilde{m}_k\} \subset
\Delta\cap \ZZ^n$.  Assume that $\Delta$ is the convex hull of
$\mathcal{A}$ and that $\mathcal{A}$ generates $\ZZ^n$ affinely over
$\ZZ$ $($meaning that $\ZZ^n$ is the $\ZZ$-span of $\{m - m' \mid m,m'
\in \mathcal{A}\})$.  Then $P_\mathcal{A} =
(x^{\tilde{m}_0},\dots,x^{\tilde{m}_k}) \in S_\Delta^{k+1}$ induces an
everywhere defined birational map
\[
\Pi : X_\Delta \to X_\mathcal{A} \subset \PP^{k}
\]
and $P_\mathcal{A}$ is the universal rational parametrization of
$X_\mathcal{A}$.
\end{proposition}

\begin{proof} 
Proposition~5.3 of \cite{sc} implies that the map $\Pi : X_\Delta \to
X_\mathcal{A}$ is the normalization map.  (In \cite{sc},
Proposition~5.3 does not assume that $\mathcal{A}$ generates $\ZZ^n$
affinely, but this is necessary since the proposition depends on
Proposition~5.2, which does assume that $\mathcal{A}$ generates the
lattice affinely.)  It follows immediately that $\Pi$ is a birational
morphism.  One can also show that $\Pi$ is strictly defined in this
case.  Then the final assertion follows immediately from
Theorem~\ref{upthm}.  This completes the proof.
\end{proof}

Here is an example of this proposition.

\begin{example}{\rm 
In the situation of Example~\ref{p1p1b2}, let $\mathcal{A} \subset
\Delta\cap \ZZ^2$ be the five vertices of $\Delta$.  Since
$\mathcal{A}$ satisfies all of the conditions of
Proposition~\ref{lattice}, it follows that
\[
P_\mathcal{A} = (x_2 x_3^2 x_4^2,  x_1^2 x_2^3 x_3^2, x_3 x_4^2 x_5, 
 x_1 x_4 x_5^2, x_1^2 x_2 x_5^2) \in S_\Delta^5
\]
is the universal rational parametrization of $X_\mathcal{A} \subset
\PP^4$.  Also note that $X_\mathcal{A}$ is the image of $X_\Delta
\subset \PP^7$ under the projection $\PP^7 -\!-\!\hskip-1.5pt \to
\PP^4$ determined by $\mathcal{A}$.}\qed
\end{example}

One final comment about Propostion~\ref{lattice} is that
$X_\mathcal{A}$ is itself a toric variety (possibly non-normal).  In
contrast, the image of $X_\Delta$ under a generic projection may fail
to be a toric variety.

\subsection{Further Examples}

Here are two examples which show what happens when we violate one of the
hypotheses of Theorem~\ref{upthm}.

\begin{example}
\label{hirzebruch}{\rm
Consider the quadrilateral $\Delta$ with vertices $(1,0)$, $(0,1)$,
$(-1,1)$ and $(-1,0)$:
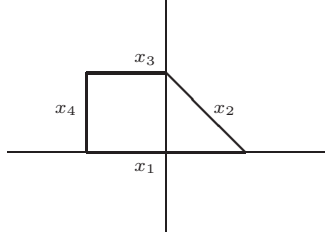
\[
\begin{matrix}
\begin{picture}(120,90)
\put(0,30){\line(1,0){120}}
\put(60,0){\line(0,1){90}}
\thicklines
\put(30,30){\line(0,1){30}}
\put(30,60){\line(1,0){30}}
\put(60,60){\line(1,-1){30}}
\put(30,30){\line(1,0){60}}
\put(48,23){$\scriptstyle{x_1}$}
\put(48,64){$\scriptstyle{x_3}$}
\put(18,45){$\scriptstyle{x_4}$}
\put(78,45){$\scriptstyle{x_2}$}
\end{picture}
\end{matrix}
\]
The lattice points $\Delta\cap \ZZ^2$ give the five monomials
\[
\begin{array}{ccc}
x_1x_2 & x_1x_4 & \\
x_2^2x_3 & x_2x_3x_4 & x_3x_4^2
\end{array}
\]
which in turn give an embedding $X_\Delta \subset \PP^4$.

Let $\mathcal{A} = \{(-1,1),(-1,0),(0,0)\} \subset \Delta\cap
\ZZ^2$.  This gives 
\[
P_\mathcal{A} = (x_1x_2, x_2^2x_3, x_2x_3x_4)
\]
and the rational map
\[
\Pi : X_\Delta -\!-\!\hskip-1.5pt \to X_\mathcal{A} \subset \PP^2
\]
defined by
\begin{equation}
\label{projdef}
(x_1x_2, x_1x_4, x_2^2x_3, x_2x_3x_4, x_3x_4^2) \mapsto (x_1x_2,
x_2^2x_3, x_2x_3x_4). 
\end{equation}
The center of this projection is the line $L = \{(0,u,0,0,v) \mid (u,v)
\ne (0,0)\}$.  One can check that $L$ is entirely
contained in $X_\Delta$ and corresponds to those points where $x_2
=0$.  Thus $\Pi$ is not strictly defined.  The surprise is that $\Pi$
is nevertheless defined everywhere on $X_\Delta$.  At first glance,
this seems impossible, since $x_2 = 0$ corresponds to points
\[
(0,x_1x_4, 0,0, x_3x_4^2) \in X_\Delta
\]
which project to $(0,0,0)$.  We get around this difficulty by letting
$x_2 = \varepsilon$, where $\varepsilon \in \CC$ is nonzero but close
to zero.  Then \eqref{projdef} becomes
\[
(x_1\varepsilon, x_1x_4, \varepsilon^2x_3, \varepsilon x_3x_4, x_3x_4^2)
\mapsto
(x_1\varepsilon, \varepsilon^2x_3, \varepsilon x_3x_4) = (x_1, \varepsilon x_3,
x_3x_4) \in \PP^2
\]
since $\varepsilon \ne 0$.  Letting $\varepsilon \to 0$, this suggests
that
\[
\Pi(0,x_1x_4, 0,0, x_3x_4^2) = (x_1, 0, x_3x_4) \in \PP^2.
\]
In fact, one can prove rigorously that $\Pi$ is defined everywhere on
$X_\Delta$.

We also note that $X_\mathcal{A} = \PP^2$ and that $\Pi$ is
birational.  This follows from
\[
\Pi^{-1}(u_0,u_1,u_2) = (u_0 u_1,u_0 u_2,u_1^2,u_1 u_2,u_2^2),
\]
where $u_0,u_1,u_2$ are homogeneous coordinates on $\PP^2$.  This is
the inverse of $\Pi$ on the open subset of $\PP^2$ where $u_0 u_1 u_2
\ne 0$.  So $\Pi$ is defined everywhere and is birational.  However,
Theorem~\ref{upthm} fails in this case.  For example, $P_\mathcal{A}$
is not a rational parametrization since $x_2$ divides the polynomials
of $P_\mathcal{A}$.  Yet the definition of rational parametrization
requires relatively prime polynomials.  Hence $P_\mathcal{A}$ has no
chance of being a universal rational parametrization.}\qed
\end{example}

Our final example concerns a singular polygon.

\begin{example}
\label{singex1}
{\rm Consider the triangle $\Delta$ with vertices
$(1,0)$, $(0,1)$, and $(-1,0)$:
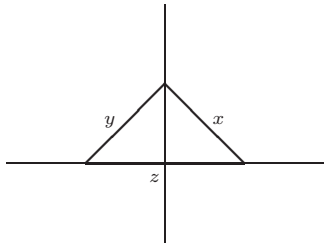
\[
\begin{matrix}
\begin{picture}(120,90)
\put(0,30){\line(1,0){120}}
\put(60,0){\line(0,1){90}}
\thicklines
\put(30,30){\line(1,1){30}}
\put(60,60){\line(1,-1){30}}
\put(30,30){\line(1,0){60}}
\put(54,23){$\scriptstyle{z}$}
\put(37,45){$\scriptstyle{y}$}
\put(78,45){$\scriptstyle{x}$}
\end{picture}
\end{matrix}
\]
The lattice points $\Delta\cap \ZZ^2$ give the four monomials
\[
\begin{array}{ccc}
& z & \\
x^2 & xy & y^2
\end{array}
\]
which in turn give an embedding $X_\Delta \subset \PP^3$ as the
singular quadric surface $u_2^2 = u_1 u_3$.  One can show that
$X_\Delta$ is the weighted projective space $\PP(1,1,2)$.

Even though $X_\Delta$ is singular, it is easy to see that
\[
P_\Delta = (z,x^2,xy,y^2)
\]
satisfies the other hypotheses of Theorem~\ref{upthm}.  So how close
is $P_\Delta$ to being a universal rational parametrization?  

For an example of what can go wrong, consider $H = (v,u,u,u)$.  This
is a rational parametrization of $X_\Delta$, yet $H$ is not of
the form $P_\Delta\circ F$ for any $F = (f_1,f_2,f_3) \in \CC[u,v]^3$
since $u$ is not a square.  So Theorem~\ref{upthm} fails in this case.

However, it is true that $H = P_\Delta\circ \widetilde{F}$, where
\begin{equation}
\label{crazy2}
\widetilde{F} = (\sqrt{u},\sqrt{u},v).
\end{equation}
So it may
be that for singular toric varieties, square roots and other radicals
appear naturally in considering what a universal parametrization
means.  But in the next section, we will learn a better method which
uses resolution of singularities.}\qed
\end{example}

\section{Universal Rational Parametrizations (Singular
Case)} 
\label{singular}

So far, we have always assumed that $X_\Delta$ is smooth, and we saw
in Example~\ref{singex1} how things can go wrong when $X_\Delta$ has
singularities.  We will use a toric resolution of singularities to
show that we still have universal rational parametrizations, where
$\Delta$ is now allowed to be \emph{any} $n$-dimensional lattice
polytope in $\ZZ^n$.

As in Section~\ref{projective}, assume that we have 
\begin{equation}
\label{pdelta}
P = (p_{0},\dots,p_{s}) \in S_\Delta^{s+1},
\end{equation}
which induces a strictly defined birational map
\[
\Pi : X_{\Delta} \to X \subset \PP^s.
\]
Our goal is to describe a universal rational parametrization of $X
\subset \PP^s$.  Our main tool will be a toric resolution of
singularities.  As shown in \cite{Fulton}, the normal fan
$\Sigma_\Delta$ of $\Delta$ has a refinement $\Sigma$ such that
$X_\Sigma$ is smooth.  It follows that the induced toric morphism
\[
\varphi : X_\Sigma \to X_\Delta
\]
is a resolution of singularities.  We may assume that $\varphi^{-1}$
is defined on the smooth part of $X_\Delta$.

Let $x_1,\dots,x_{\tilde{r}}$ be the homogeneous coordinates of
$X_\Sigma$ and let $\tilde{n}_i$ generate the edge of $\Sigma$
corresponding to $x_i$.  Some of the $\tilde{n}_i$'s will be inner
normals to facets of $\Delta$, while others will be new vectors added
in the process of resolving singularities.  We will regard the new
$\tilde{n}_i$'s as inner normals to ``virtual facet hyperplanes'' of
$\Delta$ in the following way.

Given $\tilde{n}_i$, we know that it lies in some cone $\sigma \in
\Sigma_\Delta$.  We pick the smallest such cone.  Its generators are
facet normals of $\Delta$, and the intersection of the corresponding
facets is a face $\Delta_\sigma$ of $\Delta$.  Using the support
functions defined in \cite{Fulton}, one can prove that there is a
unique integer $\tilde{a}_i$ such that
\[
\{m \in \RR^n \mid \langle m,\tilde{n}_i\rangle + \tilde{a}_i = 0\}
\cap \Delta = \Delta_\sigma.
\]
We call $\{m \in \RR^n \mid \langle m,\tilde{n}_i\rangle + \tilde{a}_i
= 0\}$ the \emph{virtual facet hyperplane} of $\Delta$ with
$\tilde{n}_i$ as inner normal.  When $\tilde{n}_i$ is the inner normal
of a facet of $\Delta$, then one easily sees that the virtual facet
hyperplane is the facet hyperplane $\langle m,\tilde{n}_i\rangle +
\tilde{a}_i = 0$ containing the corresponding facet of $\Delta$.

Let's illustrate what this looks like in one of our previous examples.

\begin{example}
\label{singex2}{\rm Consider the triangle $\Delta$ of
Example~\ref{singex1} and let $\Sigma$ be the following refinement of
its normal fan:
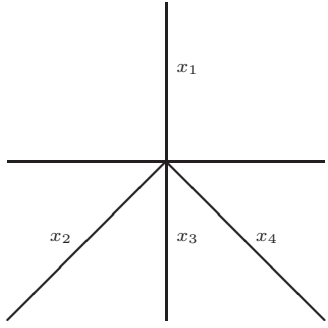
\[
\begin{matrix}
\begin{picture}(120,120)
\thinlines
\put(0,60){\line(1,0){120}}
\thicklines
\put(60,0){\line(0,1){120}}
\put(60,60){\line(-1,-1){60}}
\put(60,60){\line(1,-1){60}}
\put(94,30){$\scriptstyle{x_4}$}
\put(64,94){$\scriptstyle{x_1}$}
\put(16,30){$\scriptstyle{x_2}$}
\put(64,30){$\scriptstyle{x_3}$}
\end{picture}
\end{matrix}
\]
(So the refinement is given by adding the edge corresponding to
$x_3$.)  Let $\tilde{n}_i$ generate the edge corresponding to $x_i$.
Thus $\tilde{n}_1$, $\tilde{n}_2$ and $\tilde{n}_4$ are inner normals
of facets of the triangle $\Delta$ of Example~\ref{singex1}, while
$\tilde{n}_3$ was added to make $X_\Sigma$ smooth.  Then we can draw
the virtual facet hyperplanes (= lines in this case) and their
corresponding variables as follows:
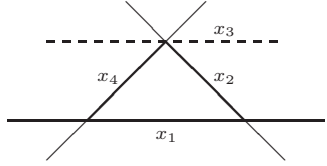
\[
\begin{matrix}
\begin{picture}(120,90)
\put(0,30){\line(1,0){120}}
\put(15,15){\line(1,1){60}}
\put(105,15){\line(-1,1){60}}
\thicklines
\put(30,30){\line(1,1){30}}
\put(60,60){\line(1,-1){30}}
\put(30,30){\line(1,0){60}}
\put(56,23){$\scriptstyle{x_1}$}
\put(34,45){$\scriptstyle{x_4}$}
\put(78,45){$\scriptstyle{x_2}$}
\put(78,63){$\scriptstyle{x_3}$}
\multiput(15,60)(6,0){15}{\line(1,0){3}}
\end{picture}
\end{matrix}
\]
The facet hyperplanes are solid lines, while the one virtual facet
hyperplane is the dashed line corresponding to $x_3$.  Note also that 
\[
\tilde{a}_1 = 0\quad \text{and} \quad \tilde{a}_2 = \tilde{a}_3 =
\tilde{a}_4 = 1.
\]
We will return to this example shortly.}\qed
\end{example}

Given this setup, a lattice point $m \in \Delta \cap \ZZ^n$ gives the
monomial 
\begin{equation}
\label{xds}
x^m = \prod_{i=1}^{\tilde{r}} x_i^{\langle m,\tilde{n}_i\rangle +
\tilde{a}_i}.
\end{equation}
We call $x^m$ a $\Delta$-monomial of the toric variety $X_\Sigma$.
Note that the exponent of $x_i$ in $x^m$ measures the lattice distance
from $m$ to the corresponding virtual facet hyperplane.  Here is an
example.

\begin{example}
\label{singex3}{\rm In the situation of Example~\ref{singex2}, the
lattice points of $\Delta\cap \ZZ^2$ give the $\Delta$-monomials 
\begin{equation}
\label{singlat}
\begin{array}{ccc}
& x_1 & \\
x_2^2x_3 & x_2x_3x_4 & x_3x_4^2
\end{array}
\end{equation}
in the homogeneous coordinates of the toric variety $X_\Sigma$ which
resolves the singularities of $X_\Delta$.}\qed
\end{example}


One useful observation is that when dealing with lattice polygons, the
only places we need to add virtual facet hyperplanes are at vertices
whose adjacent inner normals do not form a basis of $\ZZ^2$ over
$\ZZ$.  Furthermore, in this situation, there is a unique minimal
resolution of singularities which can be computed
algorithmically---see \cite[Sec.\ 2.6]{Fulton}.  Thus there is an
algorithm for finding the virtual inner normals that need to be added
at these vertices.

We are almost ready to state our main result.  As above, $\Delta$ is
an $n$-dimensional lattice polytope in $\RR^n$ and $\varphi: X_\Sigma
\to X_\Delta$ is a toric resolution.  The lattice points in
$\Delta\cap\ZZ^n$ determine
\[
S_\Delta = \mathrm{Span}(x^m \mid m \in \Delta\cap\ZZ^n)
\]
where $x^m$ is now the $\Delta$-monomial \eqref{xds} in the
homogeneous coordinates $x_1,\dots,x_{\tilde{r}}$ of $X_\Sigma$.    
Now let 
\[
P = (p_0,\dots,p_s) = \Big({\textstyle \sum_m a_{0m} x^m, \dots,
\sum_m a_{sm} x^m}\Big) \in S_\Delta^{s+1}.
\]
Then $P$ induces a rational map $p : \CC^{\tilde{r}} -\!-\!\skip1.5pt \to
\PP^s$, and similar to Proposition~\ref{pfactors}, one can show that
$p$ factors as 
\[
\CC^{\tilde{r}} -\!-\!\skip1.5pt \to X_\Sigma \xrightarrow{\ \varphi\
} X_\Delta -\!-\!\skip1.5pt \to X \subset \PP^s,
\]
where as usual, $X$ is the Zariski closure of the image of $p$.  The
map from $X_\Delta$ to $X$ will be denoted $\Pi$, and as in
Theorem~\ref{upthm}, we will assume that $\Pi$ is strictly defined and
birational.  Let $U \subset X$ be the maximal open set on which the
inverse of $\Pi: X_\Delta \to X$ is defined, and then set
\[
\widetilde{U} = U \cap \big(X \setminus \Pi(\{ x \in X_\Delta \mid x\
\text{is a singular point of}\ X_\Delta\})\big).
\]
Finally, we have $G \subset (\CC^*)^{\tilde{r}}$.  Then $\mu =
(\mu_1,\dots,\mu_{\tilde{r}}) \in G$ gives $\mu_{\Delta,\Sigma} =
\prod_{i=1}^{\tilde{r}} \mu_i^{\tilde{a}_i}$.  Let $G_{\Delta,\Sigma}$
be the kernel of the homomorphism $\mu \mapsto \mu_{\Delta,\Sigma}$.
We use this notation because $\mu_{\Delta,\Sigma}$ and
$G_{\Delta,\Sigma}$ depend not only on the polytope $\Delta$ but also
on the fan $\Sigma$.

We now show that $P$ is a universal rational parametrization of $X$.

\goodbreak

\begin{theorem}
\label{singupthm}
Let $\Delta$, $\Sigma$, ${P} = (p_0,\dots,p_s)$, $X$ and
$\widetilde{U}$ be as above.  Then $P$ is a {\bfseries universal
rational parametrization} of $X \subset \PP^s$ in the following sense:
\begin{enumerate}
\item If $F = (f_{1},\dots,f_{\tilde{r}}) \in R^{\tilde{r}} =
\CC[y_1,\dots,y_d]^{\tilde{r}}$ is $\Sigma$-irreducible, then
\[
{P}\circ F = ({\textstyle \sum_{m} a_{0m} f^{m}, \dots,\sum_{m} 
a_{sm}f^{m}}) \in R^{s+1}
\]
is a rational parametrization of $X \subset \PP^{s}$.
\item Conversely, given any rational parametrization $H \in R^{s+1}$
of $X$ whose image meets the open set $\widetilde{U} \subset
X$, there is a $\Sigma$-irreducible $F = (f_{1},\dots,f_{\tilde{r}})
\in R^{\tilde{r}}$ such that $H = {P}\circ F$.
\item If $F$ and $F'$ are $\Sigma$-irreducible, then $P\circ
F = P\circ F'$ as rational parametrizations if and only if $F' =
\mu\cdot F$ for some $\mu \in G_{\Delta,\Sigma}$.
\end{enumerate}
\end{theorem} 

The proof will be given in Section~\ref{theory}.  Note that the
theorem uses the concept of $\Sigma$-irreducible.  This uses the
obvious modification of Definition~\ref{irreddef} which applies to any
fan $\Sigma$.

Let's apply Theorem~\ref{singupthm} to the singular example we've been
studying.

\begin{example}
\label{singex4}{\rm Let $\Delta$ be the triangle of
Examples~\ref{singex1}, \ref{singex2} and~\ref{singex3}.  This gives
the singular toric variety $X_\Delta$.  The fan $\Sigma$ from
Example~\ref{singex2} gives a resolution of singularites, and the
$\Delta$-monomials $x^m$ for $m \in \Delta\cap \ZZ^2$ are given in
\eqref{singlat}.  Let
\[
P = (x_1, x_2^2x_3, x_2x_3x_4, x_3x_4^2).
\]
Since $\Pi : X_\Delta \to X \subset \PP^3$ is an isomorphism in this
case, Theorem~\ref{singupthm} implies that $P$ is a universal rational
parametrization of $X$.  

This means the following.  If $u_0,u_1,u_2,u_3$ are coordinates on
$\PP^3$, then $X$ is defined by $u_2^2 = u_1u_3$.  Hence, if $H =
(h_0,h_1,h_2,h_3)$ is a rational parametrization whose image is not
the singular point $(1,0,0,0) \in X$, then there is $F =
(f_1,f_2,f_3,f_4)$ such that
\[
H = P\circ F = (f_1,f_2^2f_3, f_2f_3f_4, f_3f_4^2).
\]
Furthermore, one can show that
\begin{itemize} 
\item $F$ is $\Sigma$-irreducible if and only if $\gcd(f_1,f_3) =
\gcd(f_2,f_4) = 1$. 
\item $F = (f_1,f_2,f_3,f_4)$ is unique up to $(f_1,\lambda
f_2,\lambda^{-2} f_3,\lambda f_4)$ for $\lambda \in \CC^*$.
\end{itemize}
For $H = (v,u,u,u) \in \CC[u,v]^4$ as in Example~\ref{singex1}, one
easily sees that $H = P\circ F$ for $F = (v,1,u,1)$ in this case.  So
unlike Example~\ref{singex1}, we don't need square roots.}\qed
\end{example}

In the smooth case, we analyzed $P$ in terms of the embedding given by
$P_\Delta$ followed by a projection. In the singular case, the analog
of $P_\Delta$ need not give an embedding.  However, when $\Delta$ is a
toric surface, then it is.  Hence the discussion of embeddings and
projections given in Section~\ref{projective} applies to \emph{any}
toric surface.

Finally, we remark that while toric resolutions are in general not
unique, in the surface case one can always find a minimal resolution
which is unique up to isomorphism.  It follows that we have a
canonical choice of universal rational parametrization when $\Delta$
is a lattice polygon.

\section{Rational Maps to Smooth Toric Varieties}
\label{rational}

In order to prove the results of Sections~\ref{projective} and
\ref{singular}, we need to study rational maps to an abstract toric
variety.  So in this section we will assume that $X_\Sigma$ is a
compact toric variety, possibly non-projective.  

In algebraic geometry, there is a well-defined notion of a rational
map between irreducible varieties, regardless of whether they are
affine, projective or defined abstractly like $X_\Sigma$.  Our goal
here is to describe \emph{all}
rational maps
\[
\CC^d -\! -\!\hskip-1.5pt \to X_\Sigma
\]
when $X_\Sigma$ is a smooth toric variety.  Recall that this means
that the generators of every $n$-dimension cone of $\Sigma$ are a
$\ZZ$-basis of $\ZZ^n$.

The natural candidate for the universal rational map to $X_\Sigma$ is
the  rational map
\[
\pi: \CC^r -\! -\!\hskip-1.5pt \to X_\Sigma
\]
of \eqref{crxs}.  So we need to explain what universal means in this
context.

Given a polynomial map $F : \CC^d \to \CC^r$ such that $F$ is
$\Sigma$-irreducible, we will show in Section~\ref{theory} that the
composition
\begin{equation}
\label{picircf}
\pi\circ F :\CC^d -\! -\!\hskip-1.5pt \to X_\Sigma.
\end{equation}
is a well-defined rational map.  One of the key assertions of
Theorem~\ref{ratmap} below is that this gives \emph{all} rational maps
from $\CC^d$ to $X_\Sigma$.

However, the map $F$ in \eqref{picircf} is not unique.  Recall from
\eqref{gdef} that we have the subgroup $G \subset (\CC^*)^r$ which is
used in the quotient representation of $X_\Sigma$.  If $F =
(f_1,\dots,f_r) \in R^r$ is $\Sigma$-irreducible and $\mu =
(\mu_1,\dots,\mu_r)\in G$, then
\[
\mu\cdot F = (\mu_1 f_1,\dots,\mu_r f_r)
\]
is also $\Sigma$-irreducible and gives the same rational map as $F$
when composed with $\pi$ (because of the quotient \eqref{quotient}).
Another key assertion of Theorem~\ref{ratmap} is that this is the
\emph{only} way that two $\Sigma$-irreducible $F$'s can give the same
$\pi\circ F$.  Thus we have complete control of the lack of
uniqueness.

We can now state the main result of this section.  Let $R =
\CC[y_1,\dots,y_d]$. 

\begin{theorem} 
\label{ratmap}
Let $X_\Sigma$ be a smooth compact toric variety.  Then:
\begin{enumerate}
\item If $F = (f_1,\dots,f_r) \in R^r$ is $\Sigma$-irreducible, then
$\pi\circ F$ gives a well-defined rational map $\pi\circ F : \CC^d -\!
-\!\hskip-1.5pt \to X_\Sigma$.
\item If $F$ and $F'$ are $\Sigma$-irreducible, then $\pi\circ F =
\pi\circ F'$ as rational maps if and only if $F' = \mu\cdot F$ for
some $\mu \in G$.
\item Finally, every rational map $f : \CC^d -\! -\!\hskip-1.5pt \to
X_\Sigma$ is of the form $\pi\circ F$ for some $\Sigma$-irreducible $F
\in R^r$.
\end{enumerate}
Hence rational maps $f : \CC^d -\!  -\!\hskip-1.5pt \to X_\Sigma$
correspond bijectively to $G$-equivalence classes of
$\Sigma$-irreducible $(f_1,\dots,f_r) \in R^r$.
\end{theorem}

The proof will be given in Section \ref{theory}.  Here is an example
of Theorem~\ref{ratmap}.

\begin{example}
\label{p1p1b}
{\rm Let $X_\Sigma$ be the toric variety of Example~\ref{p1p1b2}.
There, we saw that $F = (f_{1},\dots,f_{5})$ is $\Sigma$-irreducible
provided
\[
\gcd(f_{1},f_{3}) = \gcd(f_{1},f_{4}) = \gcd(f_{2},f_{4}) = 
\gcd(f_{2},f_{5}) = \gcd(f_{3},f_{5}) = 1
\]
and that
\[
G = \{(\lambda,\mu,\nu,\lambda \mu, \mu\nu) \mid \lambda,\mu,\nu \in 
\CC^{*}\}.
\]
By Theorem~\ref{ratmap}, it follows that rational maps to $X_{\Sigma}$
are all of the form $\pi\circ F$, where $F$ is unique up to
$(\lambda,\mu,\nu,\lambda \mu, \mu\nu)\cdot F$.

Let's look at the specific example of the map $F' : \CC^{2} \to
\CC^{5}$ defined by
\begin{equation}
\label{badmap}
F'(u,v) = (uv,1,u,v,1)
\end{equation}
This induces a rational map $\pi\circ F' : \CC^2 -\!  -\!\hskip-1.5pt
\to X_\Sigma$.  However, $F'$ is not $\Sigma$-irreducible.  To get a
$\Sigma$-irreducible representation, observe that
\begin{equation}
\label{bad2}
F' = (uv,1,u,v,1) = (uv,u^{-1},u,v,1)\cdot (1,u,1,1,1).
\end{equation}
Since $(uv,u^{-1},u,v,1) \in G$ for $u,v \ne 0$, we see that $F'$ and
$F = (1,u,1,1,1)$ give the same rational map.  Since $F$ is
$\Sigma$-irreducible, this is the representation given by
Theorem~\ref{ratmap}.  

Notice that even though \eqref{badmap} is given by polynomials,
\eqref{bad2} shows that it is not a polynomial multiple of the
$\Sigma$-irreducible representation $F = (1,u,1,1,1)$.

We can also look at \eqref{badmap} from the point of view of
Theorem~\ref{upthm} and Corollary~\ref{univcor}.  If we compose
\eqref{badmap} with $P_\Delta$ from Example~\ref{p1p1b2}, we obtain
\[
H' = (u^2 v^2, u^3 v^2, u^4 v^2, u v^2, u^2 v^2,u^3 v^2, u v^2, u^2 v^2),
\]
which satisfies the hypothesis of Corollary~\ref{univcor}.  Factoring
out the gcd $uv^2$, we can write this as
\[
H' = uv^2 (u, u^2, u^3, 1, u, u^2, 1, u) = uv^2 H.
\]
Furthermore, one easily computes that $H = P_\Delta(1,u,1,1,1)$.  Thus 
\[
H' = uv^2 P_\Delta(1,u,1,1,1).
\]
Notice also that unlike \eqref{bad2}, the representation given by
Corollary~\ref{univcor} involves only polynomials.\qed}
\end{example}  

Finally, we observe that the smoothness assumption in
Theorem~\ref{ratmap} is necessary, as shown by the following example.

\begin{example} 
\label{crazyex}
{\rm Consider the weighted projective plane $\PP(1,1,2)$.  Here,
\eqref{quotient} represents this toric variety as the quotient of
$\CC^3\setminus\{0\}$ under the action of $\CC^*$ given by
$\lambda\cdot (x,y,z) = (\lambda x,\lambda y,\lambda^2 z)$.  Then
consider the rational map
\[
\CC^2 -\! -\!\hskip-1.5pt \rightarrow \PP(1,1,2)
\]
defined by
\begin{equation}
\label{crazy}
(u,v) \mapsto (\sqrt{u},\sqrt{u},v).
\end{equation}
This looks crazy, but notice that
\[
(-1)\cdot (\sqrt{u},\sqrt{u},v) = ((-1)\sqrt{u},(-1)\sqrt{u},(-1)^2v)
= (-\sqrt{u},-\sqrt{u},v).
\]
In fact, one can prove that \eqref{crazy} gives a well-defined
rational map whose image is the curve $x = y$ in $\PP(1,1,2)$.  Since
this map cannot be written in the form $\pi\circ F$ where $F$ consists
of polynomials, we see that Theorem~\ref{ratmap} fails in this
case.

We should also mention that this example is a version of
Examples~\ref{singex1} and \ref{singex4} in disguise.  In fact, the
``crazy'' rational map \eqref{crazy} is exactly the map we used in
Example~\ref{singex1} to show that Theorem~\ref{upthm} fails for
singular toric varieties.  Recall also that we gave a purely
polynomial version of this map by using the toric resolution described
in Example~\ref{singex4}.\qed}
\end{example}

\section{Theoretical Justification}
\label{theory}

The purpose of this section is to prove the three main results of 
this paper, Theorems~\ref{singupthm}, \ref{upthm} and \ref{ratmap}.
We begin with Theorem~\ref{ratmap} since it will be used to prove the
other two theorems.

\begin{proof}[Proof of Theorem~\ref{ratmap}]
Suppose that $F = (f_1,\dots,f_r)$ is $\Sigma$-irreducible.  The
discussion following Definition~\ref{irreddef} implies that there is
some $\sigma \in \Sigma$ such that $f_{i} \ne 0$ for all $\rho_{i}
\not\subset \sigma$.  Thus $x^{\hat\sigma} = \Pi_{\rho_{i} \not\subset
\sigma} x_{i}$ does not vanish on the image of $F$, so that the image
of $F$ is not contained in the exceptional set $Z =
\mathbf{V}(x^{\hat\sigma} \mid \sigma \in \Sigma)$.  This shows that
$V = F^{-1}(Z)$ is a proper subvariety of $\CC^{d}$.  It follows that
$F$ induces a map
\[
\CC^{d} \setminus V \xrightarrow{\ F\ } \CC^r \setminus Z 
\xrightarrow{\ \pi\ } (\CC^r \setminus Z)/G = X_\Sigma.
\]
Since $\CC^{d} \setminus V$ is a nonempty Zariski open subset of
$\CC^{d}$, we get a well-defined rational map $\pi\circ F : \CC^d -\!
-\!\hskip-1.5pt \to X_\Sigma$.  This proves assertion (1) of the
theorem.

Before proving (2), we need to describe $Z$.  Since
$\mathbf{V}(x^{\hat\sigma})$ is a union of coordinate hyperplanes,
their intersection $Z$ is a union of coordinate subspaces $W_1\cup
\dots \cup W_s$.  Let one of these be $W_j =
\mathbf{V}(x_{i_1},\dots,x_{i_k})$.  Suppose that
$\rho_{i_1},\dots,\rho_{i_k} \subset \sigma$ for some $\sigma \in
\Sigma$.  Then the point $a \in \CC^r$ with $a_{i_1} = \dots = a_{i_k}
= 0$ and $a_i = 1$ otherwise lies in
$\mathbf{V}(x_{i_1},\dots,x_{i_k}) \subset Z$, yet $x^{\hat\sigma}$ is
nonvanishing at $a$.  This contradiction shows that no cone of
$\Sigma$ contains $\rho_{i_1},\dots,\rho_{i_k}$.

Note that $V = F^{-1}(Z) = F^{-1}(W_1)\cup \dots \cup
F^{-1}(W_s)$, and if we write $W_j =
\mathbf{V}(x_{i_1},\dots,x_{i_k})$ as above, then
\[
F^{-1}(W_j) = \mathbf{V}(f_{i_1},\dots,f_{i_k}).
\]
The above paragraph and Definition~\ref{irreddef} imply that
$\mathrm{gcd}(f_{i_1},\dots,f_{i_k}) = 1$.  This shows that each
$F^{-1}(W_j)$ has codimension at least 2 in $\CC^{d}$, so that the
same is true for their union $V$. 

Now we can prove uniqueness.  Suppose that another
$\Sigma$-irreducible $F' = (f'_{1},\dots,f'_{r})$ gives the same
rational map $f$.  This means the following.  Let $V' =
(F')^{-1}(Z)$.  Then as above $V'$ has codimension at least
2 and the induced rational map $f'$ is defined on $\CC^{d}\setminus
V'$.  Then $f = f'$ as rational maps implies that $f = f'$ as functions
on $U = \CC^{d}\setminus(V\cup V')$.  Since $X_\Sigma =
\big(\CC^r\setminus Z\big)/G$, this means that for each $y
\in U$, there is $\mu(y) \in G$ such that $\mu(y) \cdot F(y) = F'(y)$. 
Since $X_\Sigma$ is smooth, the quotient map
$\CC^r\setminus Z \to X_\Sigma$ is a smooth fibration with
fibers isomorphic to $G$.  This implies that the map $y \mapsto
\mu(y)$ is an algebraic map $U \to G$.

Using $G \subset (\CC^*)^r$, $\mu$ gives maps $\mu_i : U \to \CC^*$
such that $f_i(y) = \mu_i(y) f'_i(y)$ for all $y \in U$.  Now comes the
key point: since $V \cup V'$ has codimension at least 2, $\mu_i$ must
be constant.  (To see this, write $\mu_i = A/B$, where $A,B$ are
relatively prime polynomials.  Then $A$ nonconstant $\Rightarrow$ the
zeros of $\mu_i$ have codimension 1 and $B$ nonconstant $\Rightarrow$
the poles of $\mu_i$ have codimension 1.  But $\mu_i$ is defined and
nonzero outside a set of codimension at least 2.)

Hence the $\mu_i$ are constant.  It follows $\mu \in G$ is also
constant, and then $\mu$ gives the desired equivalence between
$(f_{1},\dots,f_{r})$ and $(f'_{1},\dots,f'_{r})$.  This completes the
proof of (2).

It remains to show that all rational maps from $\CC^d$ to $X_\Sigma$
arise this way.  So let $f : \CC^d -\!  -\!\hskip-1.5pt \to X_\Sigma$
be a rational map.  This is defined on some nonempty Zariski open
subset of $\CC^{d}$, and the union $U$ of all such subsets is the
maximal Zariski open subset on which $f$ is defined.  The base points
of $f$ are the complement $\CC^d \setminus U$, and since we are
mapping into a compact space, the base points have codimension at
least 2.

First consider the case when $f(U) \cap (\CC^*)^n \ne \emptyset$.
Recall that $X$ has divisors $D_1,\dots,D_r$ corresponding to
$x_1,\dots,x_r$.  Since $X_\Sigma \setminus (\CC^*)^n = \bigcup_i
D_i$, it follows that $f(U) \not\subset D_i$ for all $i$.  Thus
$f^{-1}(D_i) \subset U$ is a proper subvariety (possibly empty) for
each $i$.  If $f^{-1}(D_i) = \emptyset$, set $f_i = 1$.  Now suppose
that $f^{-1}(D_i) \ne \emptyset$.  Since $X_\Sigma$ is smooth, $D_i
\subset X_\Sigma$ is locally defined by a single equation, say $h =
0$, and then $f^{-1}(D_i) \subset U$ is defined locally by $h\circ f =
0$.  It follows that every irreducible component of $f^{-1}(D_i)$ in
$U$ has codimension 1, although the components may have
multiplicities.  Now, using $U \subset \CC^d$, we get the Zariski
closure $Z_i \subset \CC^d$ of $f^{-1}(D_i)\subset U$.  The
irreducible components of $Z_i$ also have codimension~1, with the same
multiplicities.  It follows that there is $f_i \in R$, unique up to a
constant, such that $\mathbf{V}(f_i) = Z_i$ with the same
multiplicities.

We claim that $(f_1,\dots,f_r)$ is $\Sigma$-irreducible.  Suppose that
$\rho_{i_1},\dots,\rho_{i_k}$ are contained in no cone of $\Sigma$. 
Then the relation between cones and divisors implies that $D_{i_1}
\cap \dots \cap D_{i_k} = \emptyset$ in $X_\Sigma$.  Thus, in $U$, we
have
\[
f^{-1}(D_{i_1}) \cap \dots \cap f^{-1}(D_{i_k}) = \emptyset. 
\]
Since $\mathbf{V}(f_i) \cap U = f^{-1}(D_{i})$ for all $i$, it follows
that
\begin{equation}
\label{capu}
\mathbf{V}(f_{i_1},\dots,f_{i_k}) \cap U = \emptyset.
\end{equation}
Hence $\mathbf{V}(f_{i_1},\dots,f_{i_k}) \subset \CC^d\setminus U$.
Since $\CC^d\setminus U$ has codimension at least 2, this implies that
$\mathrm{gcd}(f_{i_1},\dots,f_{i_k}) = 1$.  Thus $(f_1,\dots,f_r)$ is
$\Sigma$-irreducible.  It follows that $(c_1 f_1,\dots,c_r f_r)$ is
$\Sigma$-irreducible whenever $c_i \in \CC^*$.  This will be useful
below. 

We next claim that there are $c_i \in \CC^*$ such that $(c_1
f_1,\dots,c_r f_r)$ gives the rational map $f$.  Let $f'$ be the
rational map determined by $(f_1,\dots,f_r)$.  Using \eqref{capu} and
our earlier description of $F^{-1}(Z)$, one easily shows
that $f'$ is defined on $U$.  Furthermore, the $f_i$ were defined so
that in $U$, we have
\begin{equation}
\label{ffp}
(f')^{-1}(D_i) = f^{-1}(D_i)
\end{equation}
for all $i$.  This equality also gives the correct multiplicities.

Now take a $n$-dimensional cone $\sigma \in \Sigma$.  This gives the
affine toric variety $U_\sigma \subset X_\Sigma$, and one easily sees
that $U_\sigma = X_\Sigma \setminus \bigcup_{\rho_i \not\subset
\sigma} D_i$.  Then \eqref{ffp} implies that $(f')^{-1}(U_\sigma) =
f^{-1}(U_\sigma)$.  Call this $U'_\sigma$ and note that $U'_\sigma \ne
\emptyset$ since $f(U) \cap (\CC^*)^n \ne \emptyset$.  Thus $f$ and
$f'$ give maps $U'_\sigma \to U_\sigma$.  But since $X_\Sigma$ is
smooth, we have $U_\sigma \simeq \CC^n$.  We may assume that the edges
of $\Sigma$ are labeled so that $\rho_1,\dots,\rho_n$ are the edges of
$\sigma$.  Then write
\begin{align*}
	f\res{U'_\sigma} &= (h_1,\dots,h_n) : U'_\sigma \to \CC^n\\
	f'\res{U'_\sigma} &= (h_1',\dots,h_n') : U'_\sigma \to \CC^n.
\end{align*}
We have set things up so that $D_i\cap U_\sigma$ is defined by the
vanishing of the $i$th coordinate.  Since \eqref{ffp} respects
multiplicities, we see that $h_i'/h_i = \beta_i$ is a nonvanishing
function on $U'_\sigma$.  Thus
\[
\beta_\sigma = (\beta_1,\dots,\beta_n) : U'_\sigma \longrightarrow
(\CC^*)^n 
\]
is an algebraic map which satisfies
\[
\beta_\sigma(y) \cdot f(y) = f'(y)
\]
for all $y \in U'_\sigma$.  If $\tau$ is another $n$-dimensional cone,
then we get $\beta_\tau$ defined on $U'_\tau$ with similar properties.
However, for any $y$ in the nonempty open subset $f^{-1}((\CC^*)^n)
\subset U'_\sigma \cap U'_\tau$, there is a unique element of
$(\CC^*)^n$ which takes $f(y)$ to $f'(y)$.  It follows easily that
$\beta_\sigma = \beta_\tau$ on $U'_\sigma \cap U'_\tau$.  Furthermore,
the $U'_\sigma$ cover $U$ since the $U_\sigma$ cover $X$.  It follows
that we get an algebraic map
\[
\beta : U \longrightarrow (\CC^*)^n
\]
with the property that 
\[
\beta(y) \cdot f(y) = f'(y)
\]
for all $y \in U$.  But arguing as above, we see that $\beta$ must be
constant since $\CC^d\setminus U$ has codimension at least 2.  Thus
there is $\beta \in (\CC^*)^n$ such that $\beta\cdot f(y) = f'(y)$ for
all  $y \in U$.

Since $(\CC^*)^r/G = (\CC^*)^n$, we can pick $(c_1,\dots,c_r) \in
(\CC^*)^r$ which maps to $\beta \in (\CC^*)^n$.  We conclude that
$(c_1 f_1,\dots,c_r f_r)$ is $\Sigma$-irreducible and gives $f$.

Finally, we need to discuss what happens when our rational map $f$
satisfies $f(U) \cap (\CC^*)^n = \emptyset$.  Here, the idea is that
there is a smallest torus orbit which meets $f(U)$.  The Zariski
closure of this orbit will be $D_{i_1}\cap \cdots\cap D_{i_\ell}$
where $\rho_{i_1},\dots,\rho_{i_\ell}$ are the edges of some $\sigma_0
\in \Sigma$.  Let $\mathrm{orb}(\sigma_0)$ denote this orbit.  Then
make the following changes in the above proof:
\begin{enumerate}
\item First, let $f_{i_1} = \dots = f_{i_\ell} = 0$.
\item Second, replace $(\CC^*)^n$ with $\mathrm{orb}(\sigma_0)$  
\item Third, for $\rho_i \not\subset \sigma_0$, pick $f_i$ so that
$\mathbf{V}(f_i) \cap U = f^{-1}(D_{i})$ (with the same
multiplicities).  
\item Fourth, use $n$-dimensional cones $\sigma$ which contain
$\sigma_0$ as a face.  
\end{enumerate}
With these changes, the above argument shows that $f$ comes from a
$\Sigma$-irreducible element of $R^r$.  We omit the details.
\end{proof}

\begin{remark} {\rm In the existence part of the above proof, notice
that the set $U$ was the maximal open subset of $\CC^d$ on which $f$
was defined.  Yet the $f'$ we constructed was defined on a potentially
bigger set, namely $\CC^d\setminus F^{-1}(Z)$.  Once we
prove $\beta\cdot f = f'$, it follows that $U = \CC^d\setminus
F^{-1}(Z)$.  Using this, we obtain the following
corollary.}
\end{remark}

\begin{corollary} Let $f : \CC^d -\!  -\!\hskip-1.5pt \to X_\Sigma$ be
induced by a $\Sigma$-irreducible $F = (f_1,\dots,f_r) \in R^r$.  Then
the maximal open subset of $\CC^d$ on which $f$ is defined is given by
\[
U = \CC^d\setminus F^{-1}(Z).
\]
\end{corollary}

We now turn to the proof of the existence of universal rational 
parametrizations for smooth toric projective toric varieties.

\begin{proof}[Proof of Theorem~\ref{upthm}]  To prove (1), first
assume that $F = (f_{1},\dots,f_{r})$ is
$\Sigma_{\Delta}$-irreducible.  We need to prove that the polynomials
$\sum_{m} a_{im} f^{m}$ are relatively prime.  So suppose that an
irreducible polynomial $q \in R$ divides all of them.

We will use the interpretation of $\Pi : X_\Delta \to X$ as the
composition of the embedding $X_\Delta \subset \PP^\ell$ given by
$\Pi_\Delta$ followed by a projection.  In particular, if $L$ is the
center of the projection, then $X_\Delta \cap L = \emptyset$ since
$\Pi$ is strictly defined on $X_\Delta$.

Suppose that we have $a \in \CC^{d}$ such that $q(a) = 0$.  If $F(a)
\in \CC^{r} \setminus Z$, then $p(F(a))$ gives a point in
$X_{\Delta}\cap L$, which is empty by assumption.  It follows that
$F(\mathbf{V}(q)) \subset Z$.  Since $q$ is irreducible,
$F(\mathbf{V}(q))$ must lie in some irreducible component of $Z$.  By
the proof of Theorem~\ref{ratmap}, it follows that $F(\mathbf{V}(q))
\subset \mathbf{V}(x_{i_1},\dots,x_{i_k})$, where no cone of $\Sigma$
contains $\rho_{i_1},\dots,\rho_{i_k}$.  Thus $f_{i_{j}}$ vanishes on
$\mathbf{V}(q)$, so that $q$ divides $f_{i_{j}}$.  But this is
impossible since $F$ is $\Sigma_{\Delta}$-irreducible.  This completes
the proof of (1).

Next suppose that $H = (h_{0},\dots,h_{s})$ is a rational
parametrization of $X$ whose image meets $U \subset X$.  This gives a
rational map denoted $h :\CC^d -\!-\!\hskip-1.5pt \to X$.  Since $\Pi
: X_{\Delta} \to X$ is birational and $\Pi^{-1}$ is defined on $U$, we
get a rational map
\[
f = \Pi^{-1}\circ h : \CC^{d} -\!-\!\hskip-1.5pt \to X_{\Delta}.
\]
By Theorem~\ref{ratmap}, $f$ is induced by a
$\Sigma_\Delta$-irreducible $F = (f_{1},\dots,f_{r}) \in R^{r}$.  It
follows that $H$ and $P\circ F$ give the same rational map $\CC^{d}
-\!  -\!\hskip-1.5pt \to \PP^{s}$.  Since both satisfy the gcd
condition of Definition~\ref{rpdef}, we see that $H = c\,P\circ F$ for
some constant $c \ne 0$.

We claim that there is $\mu \in G$ such that $H = P\circ (\mu\cdot
F)$.  Recall from \eqref{pmx} that if $\mu =
(\mu_{1},\dots,\mu_{r}) \in G$, then
\begin{equation}
\label{pamud}
P(\mu\cdot (x_{1},\dots,x_{r})) = \mu_{\Delta}
P(x_{1},\dots,x_{r}), 
\end{equation}
where 
\[
\mu_{\Delta} = \prod_{i=1}^{r} \mu_{i}^{a_{i}}.
\]
Assume for the moment that the map
\begin{equation}
\label{gdelta}
G \longrightarrow \CC^{*}
\end{equation}
defined by $\mu \mapsto \mu_{\Delta}$ is surjective.  Then we can find
$\mu \in G$ such that $\mu_{\Delta} = c$.  It follows that
\[
H = c\,P\circ F = \mu_{\Delta} P\circ F =
P\circ (\mu\cdot F),
\]
as claimed.  Since $G_\Delta$ is the kernel of \eqref{gdelta}, the
uniqueness assertion of Theorem~\ref{ratmap} easily implies that
$\mu\cdot F$ is unique up to $G_\Delta$-equivalence.

We still need to prove that \eqref{gdelta} is surjective.  Since this
map is a character, its image is either finite or all of $\CC^*$.
Furthermore, it is well-known that $G$ is connected since $X_{\Delta}$
is smooth.  Hence the image is either the identity or $\CC^*$.  So all
we need to prove is that \eqref{gdelta} is nonconstant.

If the map is constant, then $\mu_{\Delta} = 1$ for all $\mu \in G$.
We claim this implies the existence of $m \in M$ such that
\begin{equation}
\label{badm}
\langle m,n_i\rangle = a_i\quad\text{for all}\ i = 1,\dots,r.
\end{equation}
We prove this as follows.  As explained in \cite{hc}, the inclusion
$G \subset (\CC^{*})^{r}$ induces an exact sequence
\[
1 \longrightarrow G \longrightarrow (\CC^{*})^{r} \xrightarrow{\ \phi\ } 
(\CC^{*})^{n} \longrightarrow 1.
\]
The map $\mu \to \mu_{\Delta} = \prod_{i=1}^{r} \mu_{i}^{a_{i}}$
extends to the character $(\CC^{*})^{r} \to \CC^{*}$ corresponding to
$(a_1,\dots,a_r) \in \ZZ^r$.  If $\mu \to \mu_{\Delta}$ is constant on
$G$, then above exact sequence shows that it induces a character
$\chi^{m} : (\CC^{*})^{n} \to \CC^{*}$.  Since the map $\phi$ is dual
to the inclusion $\ZZ^{n} \to \ZZ^{r}$ which sends $m$ to $(\langle
m,n_{1}\rangle, \dots, \langle m,n_{r}\rangle)$, it follows that
$(a_{1},\dots,a_r) = (\langle m,n_{1}\rangle,\dots,\langle
m,n_{r}\rangle)$, as claimed.  Thus \eqref{badm} is proved.

However, if we compare \eqref{badm} to \eqref{deltadesc}, we see that
$-m$ lies in every facet of $\Delta$, which is clearly impossible.
This contradiction shows that \eqref{gdelta} must be nonconstant, and
we are done.
\end{proof}

Finally, we prove the existence of universal rational parametrizations
for arbitrary projective toric varieties.

\begin{proof}[Proof of Theorem~\ref{singupthm}]
Recall that $P$ induces a rational map $p : \CC^{\tilde{r}}
-\!-\!\skip1.5pt \to \PP^s$ which factors
\[
\CC^{\tilde{r}} -\!-\!\skip1.5pt \to X_\Sigma \xrightarrow{\ \varphi\
} X_\Delta \xrightarrow{\ \Pi\ } X \subset \PP^s.
\]
Furthermore, the argument preceding the statement of
Theorem~\ref{upthm} shows that $p$ is a rational parametrization of
$X$.  From here, the proof of (1) is identical to the proof of the
first part of Theorem~\ref{upthm}.  

As for (2), observe that the composition
\[
X_\Sigma \xrightarrow{\ \varphi\ } X_\Delta \xrightarrow{\ \Pi\ } X
\]
is birational.  Furthermore, since $\Pi^{-1}$ is defined on $U$ and
$\varphi^{-1}$ is defined on the smooth part of $X_\Delta$, it follows
that 
\[
\widetilde{\Pi} : X_\Sigma \to X
\]
is a birational morphism whose inverse is defined on $\widetilde{U}$.
Hence, if a rational parametrization $H$ induces a rational map $h :
\CC^{d} -\!-\!\skip1.5pt \to X \subset \PP^s$ whose image meets
$\tilde{U}$, then
\[
f = \widetilde{\Pi}^{-1} \circ h : \CC^{d} -\!-\!\skip1.5pt \to
X_\Sigma
\]
is a rational map.  As in the proof of Theorem~\ref{upthm},
Theorem~\ref{ratmap} implies that $f$ is given by a
$\Sigma$-irreducible $F$.  This means that $H$ and $P\circ F$ agree up
to a constant.

For the third assertion of the theorem, the proof follows from what we
did in the proof of Theorem~\ref{upthm}.  This completes the proof of
the theorem.
\end{proof}


\begin{thebibliography}{99}

\bibitem{functor} D.\ Cox, \emph{The functor of a smooth toric
variety}, Tohoku Math.\ J. {\bf 47} (1995), 251--262.

\bibitem{hc} D.\ Cox, \emph{The homogeneous coordinate ring of a toric
variety}, J.\ Algebraic Geom. {\bf 4} (1995), 17--50.

\bibitem{sc} D.\ Cox, \emph{Recent developments in toric geometry}, in
{\sl Algebraic Geometry $($Santa Cruz, 1995\/$)$} (J.\ Koll\'ar, R.\
Lazarsfeld and D.\ Morrison, editors), Proc.\ Symp.\ Pure Math.\ {\bf
62.2}, AMS, Providence, RI, 1997, 389--436.

\bibitem{what} D.\ Cox, \emph{What is a toric variety?},
\url{http://www.cs.amherst.edu/~dac/lectures/tutorial.ps}

\bibitem{DHJ} R.\ Dietz, J.\ Hoschek, B.\ J\"utter, 
\emph{An algebraic approach to curves and surfaces on the sphere and
other quadrics}, 
Computer Aided Geometric Design\ {\bf 10} (1993), 211--229.

\bibitem{Fulton} W.\ Fulton, {\sl Introduction to Toric Varieties\/},
Princeton U.\ Press, Princeton, NJ, 1993.

\bibitem{Guest} M.\ A.\ Guest, \emph{The topology of the space of rational 
curves on a toric variety}, Acta Math.\ {\bf 174} (1995), 119--145.

\bibitem{Jac} K.\ Jaczewski, \emph{Generalized Euler sequence and
toric varieties}, in {\sl Classification of Algebraic Varieties\/}
(C.\ Cilberto, E.\ Livorni and A.\ Sommese, editors), AMS, Providence,
RI, 1994, 227--247.

\bibitem{rimas} R.\ Krasauskas, \emph{Toric surface patches},
Adv.\ Comput.\ Math.\ {\bf 17} (2002), 89--113.

\bibitem{up} R.\ Krasauskas, \emph{Universal parametrizations of some
rational surfaces}, in {\sl Curves and Surfaces with Applications in
CAGD} (A.\ Le Mehaute, C.\ Rabut and L.\ L.\ Schumaker, editors),
Vanderbilt University Press, Nashville, 1997, 231--238.

\bibitem{M} C.\ M\"aurer, \emph{Generalized parameter representations
of tori, Dupin cyclides and supercyclides}, 
in {\sl Curves and Surfaces with Applications in CAGD} 
(A.\ Le Mehaute, C.\ Rabut and L.\ L.\ Schumaker, editors),
Vanderbilt University Press, Nashville, 1997, 295--302.

\bibitem{Mu} R.\ M\"uller, \emph{Universal parametrization and
interpolation on cubic surfaces}, Comput.\ Aided Geom.\ Design {\bf
19} (2002), 479--502.

\bibitem{zube} S.\ Zub\.e, \emph{The $n$-sided toric patches and
$\mathcal{A}$-resultants}, Comput.\ Aided Geom.\ Des.\ {\bf 17} (2000),
695--714. 
\end{thebibliography}
\end{document}